\documentclass{elsarticle}

\usepackage{lineno}
\usepackage[hidelinks]{hyperref}
\usepackage{fancyhdr}
\usepackage{ifthen}

\usepackage{amsmath}
\usepackage{amscd}
\usepackage{amsfonts}
\usepackage{amssymb}
\usepackage{amsthm}     
\usepackage{mathabx}    
\usepackage{nicefrac}   

\usepackage{enumitem}   

\usepackage{tikz}               
\usepackage{pgfplots}           
\usepackage{graphicx}           
\pgfplotsset{compat=newest}
\usetikzlibrary{patterns}

\usepackage{xcolor}             

\usepackage[font={footnotesize,it},labelfont=bf]{caption}
\usepackage{float}              

\usepackage{subfig}             

\usepackage{chngcntr}           

\modulolinenumbers[5]

\journal{Journal of \LaTeX\ Templates}

\makeatletter
\def\user@resume{resume}
\def\user@intermezzo{intermezzo}
\newcounter{previousequation}
\newcounter{lastsubequation}
\newcounter{savedparentequation}
\renewenvironment{subequations}[1][]{%
      \def\user@decides{#1}%
      \setcounter{previousequation}{\value{equation}}%
      \ifx\user@decides\user@resume 
           \setcounter{equation}{\value{savedparentequation}}%
      \else  
      \ifx\user@decides\user@intermezzo
           \refstepcounter{equation}%
      \else
           \setcounter{lastsubequation}{0}%
           \refstepcounter{equation}%
      \fi\fi
      \protected@edef\theHparentequation{%
          \@ifundefined {theHequation}\theequation \theHequation}%
      \protected@edef\theparentequation{\theequation}%
      \setcounter{parentequation}{\value{equation}}%
      \ifx\user@decides\user@resume 
           \setcounter{equation}{\value{lastsubequation}}%
         \else
           \setcounter{equation}{0}%
      \fi
      \def\theequation  {\theparentequation  \alph{equation}}%
      \def\theHequation {\theHparentequation \alph{equation}}%
      \ignorespaces
}{%
  \ifx\user@decides\user@resume
       \setcounter{lastsubequation}{\value{equation}}%
       \setcounter{equation}{\value{previousequation}}%
  \else
  \ifx\user@decides\user@intermezzo
       \setcounter{equation}{\value{parentequation}}%
  \else
       \setcounter{lastsubequation}{\value{equation}}%
       \setcounter{savedparentequation}{\value{parentequation}}%
       \setcounter{equation}{\value{parentequation}}%
  \fi\fi
  \ignorespacesafterend
}

\newtheorem{theorem}{Theorem}
\theoremstyle{plain}

\newtheorem{proposition}{Proposition}

\numberwithin{equation}{section}

\usetikzlibrary{plotmarks} 
\usetikzlibrary{external}
\tikzsetexternalprefix{/}
\pgfplotsset{compat=newest}
\usetikzlibrary{patterns}
\pgfplotsset{plot coordinates/math parser=false}
\newlength\figurewidth
\newlength\figureheight

\newcommand{\e}{\operatorname{e}}
\newcommand{\real}[1]{\operatorname{Re}\left\{#1\right\}}
\newcommand{\imag}[1]{\operatorname{Im}\left\{#1\right\}}

\newcommand{\ucmcg}{u_{h}}
\newcommand{\ulcmcg}{u_{h}^{(\ell)}}
\newcommand{\uddm}{u_{h}^{*}}
\newcommand{\uddmHDG}{u_{h,HDG}^{*}}

\newcommand{\ucmcgHDG}{u_{h,HDG}}

\newcommand{\ulasymp}{w_h^{(\ell)}}
\newcommand{\laplace}{\Delta}

\newcommand{\intOmega}{\int_{\Omega}\hspace*{-0.25em}}
\newcommand{\intGammaS}{\int_{\Gamma_S}\hspace*{-0.75em}}
\newcommand{\intGammaD}{\int_{\Gamma_D}\hspace*{-0.75em}}
\newcommand{\intGammaN}{\int_{\Gamma_N}\hspace*{-0.75em}}
\newcommand{\intPartialOmega}{\int_{\partial \Omega}\hspace*{-0.5em}}
\newcommand{\intT}{\int_{0}^{T}\hspace*{-0.5em}}

\begin{document}

\begin{frontmatter}

\title{Parallel Controllability Methods \\ For the Helmholtz Equation}
%
\author[Basel]{Marcus J. Grote}
\ead{marcus.grote@unibas.ch}
\author[Paris]{Fr{\'e}d{\'e}ric Nataf}
\ead{nataf@ann.jussieu.fr}
\author[Basel]{Jet Hoe Tang}
\ead{jet.tang@unibas.ch}
\author[Paris]{Pierre-Henri Tournier}
\ead{tournier@ljll.math.upmc.fr}
\address[Basel]{University of Basel, Spiegelgasse 1, 4051 Basel, Switzerland}
\address[Paris]{Laboratoire J.L. Lions, Universit{\'e} Pierre et Marie Curie, 4 place Jussieu, 75005 Paris, France, and ALPINES INRIA, Paris, France}

\begin{abstract}
The Helmholtz equation is notoriously difficult to solve with standard numerical methods, 
increasingly so, in fact, at higher frequencies. 
Controllability methods instead transform the problem back to the time-domain, where they
seek the time-harmonic solution of the corresponding time-dependent wave equation. 
Two different approaches are considered here based either on the first or
second-order formulation of the wave equation. Both are extended to general boundary-value problems
governed by the Helmholtz equation
and lead to robust and inherently parallel algorithms. Numerical results illustrate the accuracy, convergence
and strong scalability of controllability methods for the solution of high frequency Helmholtz 
equations with up to a billion unknowns on massively parallel architectures. 
\end{abstract}

\begin{keyword}
Helmholtz equation; time-harmonic scattering; 
exact controllability; 
finite elements; 
domain decomposition; parallel scalability
\end{keyword}

\end{frontmatter}



\section{Introduction}
The efficient numerical solution of the Helmholtz equation is fundamental to 
the simulation of
 time-harmonic wave phenomena in acoustics, electromagnetics or elasticity.
As the time frequency $\omega>0$ increases, so does the
size of the linear system resulting from any numerical discretization in order
to resolve the increasingly smaller wave lengths. With the increase in frequency, however,
the performance of standard 
preconditioners based on multigrid, incomplete factorization or domain decomposition 
approaches, originally developed for positive definite Laplace-like equations, rapidly
deteriorates \citep{EG2005}. 

In recent years, a growing number of increasingly sophisticated preconditioners has
been proposed for the iterative solution of the Helmholtz equation; "Shifted Laplacian" 
preconditioners~\cite{ERLANGGA2004409}, for instance, have led to 
modern multigrid~\cite{Calandra:2017:GMP,BGS2009} and domain decomposition preconditioners~\cite{graham2017domain,BDGST2017}.
While some of those
preconditioners may achieve a desirable frequency independent convergence behavior in special situations \cite{EY2011}, that optimal behavior is often lost in the presence of strong
heterogeneity. Moreover, they are typically tied to a special discretization or fail to 
achieve optimal scaling on parallel architectures.

Controllability methods (CM) offer an alternative approach for the 
numerical solution of the Helmholtz equation. Instead of solving 
the problem directly in the frequency domain, we 
first transform it back to the time domain where we seek the corresponding
time-dependent periodic solution, $y(\cdot,t)$, with known period $T = 2\pi/\omega$. 
By minimizing an energy functional $J(v_0,v_1)$ which penalizes the mismatch after 
one period, controllability methods iteratively adjust
the (unknown) initial condition $(v_0,v_1)$ thereby steering $y(\cdot,t)$ towards the desired periodic solution. 
Once the minimizer of $J$ has been found, we immediately recover from it the solution
of the Helmholtz equation. As the CM combines the numerical integration of
the time-dependent wave equation with a conjugate gradient (CG) iteration, it is remarkably 
robust and inherently parallel.

In  \cite{BGP1998}, Bristeau et al. proposed the first CM for sound-soft scattering problems
based on the wave equation in standard
second-order form. Since the initial condition $(v_0,v_1)$ then lies in $H^1 \times L^2$, the
original formulation requires the solution of a coercive elliptic problem at each CG iteration.
Heikkola et al. in \cite{HMPR2007,HMPR2007_2} presented a higher-order 
version by using spectral FE and the classical fourth-order Runge-Kutta (RK) method. 
For more general boundary-value problems, such as wave scattering from sound-hard obstacles, inclusions,
or wave propagation in physically bounded domains, the original CM will generally fail because the
minimizer of $J$ is no longer unique. In \cite{GT2018}, we proposed alternative energy functionals which 
restore uniqueness, albeit at a small extra computational cost, for general boundary-value problems governed by the Helmholtz equation.

More recently, Glowinski and Rossi  \cite{GR2006} proposed a CM based on the wave
equation in first-order (or mixed) form using classical Raviart-Thomas (RT) finite
elements. As $(v_0,v_1)$ then
lies in $L^2 \times (L^2)^d$, the solution of an elliptic problem at each CG iteration 
is no longer necessary and the CM becomes in principle trivially parallel. Still, the lack of availability
of mass-lumping for RT elements again nullifies the main advantage of the first-order formulation
because the mass-matrix now needs to be "inverted" at each time-step.

Here we revisit the original CM from \cite{BGP1998,GR2006} and consider two distinct discretizations,
which both lead to highly efficient and inherently parallel methods. In Section 2, we recall the CMCG method
based on the wave equation in second-order form and propose a filtering procedure which permits the use of the original energy functional $J$, regardless of the boundary conditions.
Next, in Section 3, we consider the CM based on the wave equation in first-order form and again
show how to extend it to arbitrary boundary-value problems governed by the Helmholtz equation. 
Thanks to a recent hybrid discontinuous Galerkin (HDG) method \citep{CNPS2016}, 
which automatically yields a block-diagonal mass-matrix, the time integration of
the wave equation then becomes truly explicit and the entire CMCG approach trivially parallel. 
In Section 4, we perform a series of numerical experiments to illustrate the accuracy, convergence behavior 
and inherent parallelism of the CMCG approach. In particular, we apply it to large-scale high-frequency
Helmholtz problems with up to a billion unknowns to demonstrate its
strong scalability on massively parallel architectures.

\section{Controllability methods for the second-order formulation \label{sec:2nd-order-formulation}}
\subsection{Time-harmonic waves}
We consider a time-harmonic wave field $u(x)$ in a bounded connected 
computational domain $\Omega\subset \mathbb{R}^d$, $d\leq 3$, 
with a Lipschitz boundary $\Gamma$. The boundary 
consists of three disjoint components, 
$\Gamma = \Gamma_D \cup \Gamma_N \cup \Gamma_S$ where
we impose a Dirichlet, Neumann and impedance (or Sommerfeld-like absorbing) boundary condition, respectively; the boundary condition is omitted whenever the corresponding component is empty.
In $\Omega$, the wave field $u$ hence satisfies the {\it Helmholtz equation}
\begin{subequations}
  \label{eqn:HE}
  \begin{alignat}{4}
    -\laplace u(x) \ - k^2(x)\ u(x)      
  &\ = \ f(x)      ,&&\qquad  x\in\Omega,
  \label{eqn:HE_PDE}
  \\  
  \frac{\partial u(x)}{\partial n} \ -\  i k(x)\ u(x) &\ = \ g_S(x) ,&&\qquad x\in\Gamma_S
  \label{eqn:HE_ABC}
  ,  
  \\   
  \frac{\partial u(x)}{\partial n} &\ = \ g_N(x) ,&&\qquad x\in\Gamma_N
  \label{eqn:HE_NBC}
  ,
  \\
  u(x) &\ = \ g_D(x) ,&&\qquad x\in\Gamma_D
  \label{eqn:HE_DBC}
  ,
  \end{alignat}
\end{subequations}
where $\omega>0$ is the (angular) frequency,
$c(x) > 0$ the wave speed, 
$k(x) = \omega/c(x)$ the wave number,
$n$ the unit outward normal, 
and
$f$, $g_N$, $g_S$ and $g_D$ are known and may vanish.

The above formulation is rather general and encompasses most common applications such as
{\it sound-soft scattering} problems with $\Gamma_S\neq \emptyset$
and $\Gamma_D \neq \emptyset$, {\it sound-hard scattering} problems with $\Gamma_S\neq \emptyset$
and $\Gamma_N \neq \emptyset$, or {\it physically bounded domains} with $\Gamma_S=\emptyset$.
We shall always assume for any particular choice of $\omega$, $c(x)$, or combination of boundary conditions that (\ref{eqn:HE}) is well-posed and has a unique solution $u\in H^1(\Omega)$.

Instead of solving the Helmholtz equation directly in the frequency domain, we now reformulate (\ref{eqn:HE}) in the time domain.
Then, the corresponding time-harmonic wave field, $\real{u(x)\e^{-i\omega t}}$, satisfies 
the (real-valued) {\it time-dependent wave equation}
\begin{subequations}  
\label{eqn:WE}
\begin{alignat}{4} 
  \frac{1}{c^{2}(x)}\frac{\partial^2y(x,t)}{\partial^2 t} \ - \ \laplace y(x,t)
    &\ = \  \real{f(x)\e^{-i\omega t}} 
    ,&&\quad x\in\Omega, \ t>0,
  \label{eqn:WE_INT}
  \\
  \dfrac{\partial y(x,t)}{\partial n} +  \dfrac{1}{c(x)}\  \dfrac{\partial y(x,t)}{\partial t} &\ = \  \real{g_S(x)\e^{-i\omega t}} 
    ,&&\quad x\in\Gamma_S, t>0,
  \label{eqn:WE_ABC}
  \\
  \dfrac{\partial y(x,t)}{\partial n} &\ = \  \real{g_N(x)\e^{-i\omega t}} 
    ,&&\quad x\in\Gamma_N, t>0,
  \label{eqn:WE_NBC}
  \\
  y(x,t) &\ = \  \real{g_D(x)\e^{-i\omega t}} 
    ,&&\quad x\in\Gamma_D, t>0,
  \label{eqn:WE_DBC}
  \\
  y(x,0) \ = \ v_0(x),
  \quad
  \frac{\partial y(x,0)}{\partial t} &\ = \ v_1(x)
    ,&&\quad x\in\Omega,
   \label{eqn:WE_INIT}
\end{alignat}
\end{subequations}
for the (unknown) initial values $v_0 = \real{u}$ and  $v_1 = \omega\imag{u}$.

For sound-soft scattering problems (\ref{eqn:HE}), where $|\Gamma_D|>0$ and $|\Gamma_S|>0$,
Bristeau et al. \cite{BGP1998,L1988} proposed 
to determine $u(x)$ via controllability
by computing a time-periodic solution $y(x,t)$ of (\ref{eqn:WE}) with period $T={2\pi}/{\omega}$.
Once the initial values $v_0,v_1$ of $y$ are known, the solution $u$ of the original Helmholtz equation (\ref{eqn:HE}) is immediately given by
\begin{equation}
  \label{eqn:TH_SOLUTION}
  u = v_0 + \frac{i}{\omega} v_1, \qquad v_0,v_1 \in H^1(\Omega).
\end{equation}
To determine $v_0$ and $v_1$, the problem is
reformulated as a least-squares optimization problem 
over $H^1(\Omega) \times L^2(\Omega)$ for the quadratic cost functional
\begin{equation}
  \label{eqn:cost_functional_J}  
  J (v_0,v_1) =   
  \frac{1}{2} \intOmega |\nabla y(x,T)-\nabla v_0(x)|^2 dx
  +
  \frac{1}{2} \intOmega \frac{1}{c^{2}(x)}(y_t(x,T)-v_1(x))^2 dx
  ,
\end{equation}
where $y$ satisfies (\ref{eqn:WE}) with the initial values $v_0$ and $v_1$. 
The functional $J$ measures in the energy norm 
the mismatch between the solution of (\ref{eqn:WE}) at the initial time and after
one period. It is non-negative and convex, while
$J(v_0,v_1) = 0$ if, and only if, 
$\nabla y(\cdot,T)=\nabla y(\cdot,0)$ 
and $y_t(\cdot,T)=y_t(\cdot,0)$ for any given initial values $(v_0,v_1)$;
in particular, $J(v_0,v_1)=0$ if $v_0 = \real{u}$ and  $v_1 = \omega\imag{u}$.

For more general scattering problems, however, $J$ is no longer strictly convex as the $T$-periodicity of $y_t$ and $\nabla y$ no longer guarantees a unique periodic solution $y$ of (\ref{eqn:WE}). Instead, 
for the general boundary-value problem (\ref{eqn:HE}), the situation is more complicated and summarized in the following theorem \cite{GT2018} .

\begin{theorem}
\label{thm:TP_UNIQUE}
Let $u\in H^1(\Omega)$ be the unique solution of (\ref{eqn:HE}) and
$y\in C^0([0,T];H^1(\Omega))\cap C^1([0,T];L^2(\Omega))$
be a (real-valued) solution of (\ref{eqn:WE}) with initial values $(v_0, v_1)\in H^1(\Omega)\times  L^2(\Omega)$. If $\nabla y$ and $y_t$ are time periodic with period $T = 2\pi/\omega$, then
$y$ admits the Fourier series representation
\begin{equation}
\label{eqn:FourierSeries_Y}
    (y(\cdot,t),\varphi)
    \ =\  (\real{u \e^{-i\omega t}},\varphi) + 
            (\lambda+\eta t,\varphi)
            + \sum_{|\ell|>1}
            (\gamma_\ell\e^{i\omega \ell t},\varphi)
\end{equation}
for any $\varphi\in H_D^1$,
where the constants $\lambda,\eta\in\mathbb{R}$ 
and the eigenfunctions $\gamma_\ell = {\alpha}_\ell+i{\beta}_\ell$, ${\alpha}_\ell, {\beta}_\ell \in H^1(\Omega)$, $|\ell|>1$ satisfy
\begin{subequations}  
\label{eqn:Eigenproblem}
\begin{alignat}{4} 
    - \laplace\gamma_\ell(x)
    - (\ell k(x))^2  \gamma_\ell(x)
    &\ =\ 0, &&\qquad x\in\Omega, 
    \label{eqn:EP_INT}
    \\
     \frac{\partial \gamma_\ell(x)}{\partial n}
    +
    i\ell k(x)\ \gamma_\ell(x)
    &\ = \  0, &&\qquad x\in\Gamma_S,
  \label{eqn:EP_ABC}
  \\
    \frac{\partial \gamma_\ell(x)}{\partial n}
    &\ = \  0, &&\qquad x\in\Gamma_N,
  \label{eqn:EP_NBC}
  \\
  \gamma_\ell(x) 
    &\ = \  0, &&\qquad x\in\Gamma_D,
  \label{eqn:EP_DBC}  
\end{alignat}

\end{subequations}
Let $v=v_0+({i}/{\omega}) \ v_1$. Then $v$ satisfies
\begin{equation}
  \label{eqn:FourierSeries_V}
     (v,\varphi) \ =\  (u,\varphi) + (\lambda+\frac{i}{\omega}\eta,\varphi) + \sum_{|\ell|>1}(\alpha_\ell + i\ell\beta_\ell, \varphi),  \quad \forall \varphi \in H_D^1.
\end{equation}
Furthermore, if $|\Gamma_S|>0$, then $\eta = 0$. If $|\Gamma_D|>0$, then $\lambda=\eta=0$. \\
Here $H_D^1:=\{w\in H^1(\Omega): w=0 \text{ on $\Gamma_D$}\}$ and $(\cdot,\cdot)$ denotes the standard $L^2(\Omega)$ inner product.
\end{theorem}
\proof See \cite{GT2018}.
\qed

For sound-soft scattering problems ($|\Gamma_S|>0,|\Gamma_D|>0$), where both Dirichlet and Sommerfeld-like absorbing boundary conditions are imposed on $\Gamma$, all
the eigenfunctions $\gamma_\ell$, $|\ell|>1$, and the constants $\lambda$, $\eta$ in (\ref{eqn:FourierSeries_V})
vanish identically. Thus, the minimizer $v = v_0 + ({i}/{\omega}) v_1$ of $J$ in (\ref{eqn:cost_functional_J})
then coincides with $u$.

For scattering problems from sound-hard obstacles or penetrable inclusions ($|\Gamma_S|>0$, $|\Gamma_D|=0$),
the eigenfunctions $\gamma_\ell$ and the constant $\eta$ in (\ref{eqn:FourierSeries_V}) still vanish identically,
yet the constant $\lambda$ may be nonzero. Given any minimizer $v = u +\lambda$ of $J$, we can recover
$u$ by subtracting the spurious shift $\lambda$ using the compatibility condition:
\begin{equation*}
\label{eqn:CC}
\lambda 
=
  \frac{1}{\|k\|_{L^2(\Omega)}^2 +  i |k|_{L^1(\Gamma_S)} }
  \bigg(
     \intOmega k^2 v  + i\intGammaS k v
    +
    \intOmega f 
    + \intGammaS g_S+ \intGammaN g_N
  \bigg)  
  .
\end{equation*}
In fact, any impedance condition (\ref{eqn:HE_ABC}) that includes a positive (or negative) definite zeroth order term, such as a more accurate absorbing boundary condition \cite{BGT1982,GK1995}, also circumvents the indeterminacy due to $\lambda$.
For wave propagation in physically bounded domains ($|\Gamma_S|=0$), the eigenfunctions
$\gamma_\ell$ and the constants $\lambda,\eta$ in (\ref{eqn:FourierSeries_V}) typically do not vanish.
However, we can always restore uniqueness by replacing $J$ 
with an alternative energy functional,
thereby incurring a small increase in computational cost -- see  \cite{GT2018}.

\subsection{Fundamental frequency extraction via filtering
 \label{sec:FilterFundamentalFrequency}}
From Theorem \ref{thm:TP_UNIQUE} we conclude that a minimizer of $J$ generally yields 
a time-dependent
solution $y$ of  (\ref{eqn:WE}), which contains a constant shift determined by $\lambda$, a linearly growing
part determined by $\eta$, and higher frequency harmonics determined by $\gamma_\ell$, all superimposed on
the desired time-harmonic field $u$ with fundamental frequency $\omega$. Those spurious modes can be
eliminated by replacing $J$ with an alternative energy functional
at a small extra computational cost \cite{GT2018}. Instead we now 
propose an alternative approach via filtering which removes all spurious modes without requiring a modified energy functional. 

Let $y(x,t)$ be the time-dependent solution of (\ref{eqn:WE}) that corresponds to a minimizer $(v_0,v_1)$ of $J$. Next, we define $\widehat{y}\in\{w \in H^1(\Omega) \ | \ w = g_D \text{ on $\Gamma_D$}  \}$ as
\begin{equation}
  \label{eqn:fundamental-frequency-solution}
  \widehat{y}(x) := 
  \frac{1}{T}\int_{0}^{T}
    \big(
        y(x,t)+\frac{i}{\omega}y_t(x,t)
    \big)\e^{i\omega t}\,dt .
\end{equation}
To extract $u(x)$ from $y(x,t)$, we now take advantage of the mutual orthogonality of different time harmonics $\exp(i\omega \ell t)$ in $L^2(0,T)$. Hence, we multiply (\ref{eqn:FourierSeries_Y}) with $\e^{i\omega t}$  and 
 integrate in time over $(0,T)$ to obtain
\begin{alignat}{4}
 \widehat{y}(x)
&\ =\ 
\frac 1T \intT 
\big(
\operatorname{Re}\{u\e^{-i\omega t}\}
+
\lambda
+
\eta t
+
i
\operatorname{Im}\{u\e^{-i\omega t}\}
+\frac{i\eta}{\omega}
\big)\e^{i\omega t}
 dt
\nonumber\\
&\ =\ 
\frac 1T\int_{0}^{T}
u\e^{-i\omega t}
\e^{i\omega t}
dt
-\frac{i\eta}{\omega}
\ =\ 
u
-\frac{i\eta}{\omega}. \end{alignat}
This yields
\begin{equation}
\label{eqn:shifted-fundamental-frequency-solution}
 u(x) = \widehat y(x) + \frac{i\eta}{\omega}, \qquad x \in\Omega
\end{equation}
where $\lambda$ and all $\gamma_\ell$ have vanished but the constant $\eta$ is still undetermined. 

If $|\Gamma_S|>0$ or $|\Gamma_D|>0$,
Theorem \ref {thm:TP_UNIQUE} implies that $\eta=0$ and thus $u(x) = \widehat y(x)$. 
Otherwise in the pure Neumann case ($\Gamma = \Gamma_N$), we 
determine $\eta$  by integrating (\ref{eqn:shifted-fundamental-frequency-solution}), multiplied by $k^2(x)$,
over $\Omega$ and using the compatibility condition
\begin{alignat}{4}
    - \int_{\Omega} k^2(x) u(x) \ dx     
    =    
    \int_{\Omega} f(x)\ dx 
    + 
    \intPartialOmega g_N(x)\ ds
    \label{eqn:cc_neumann}.    
\end{alignat}
from (\ref{eqn:HE_PDE}). This immediately yields the remaining constant
\begin{alignat}{4}            
    \frac{i\eta}{\omega}
    &\ = \     
    -
    \frac{1}{\|k\|_{L^2(\Omega)}^2}
    \bigg(
    \int_{\Omega} f(x)\ dx    
    +
    \intPartialOmega g_N(x)\ ds    
    +
    \int_{\Omega}
     k^2(x)\widehat y(x) dx
    \bigg)
    .
    \label{eqn:cc_neumann_eta}    
\end{alignat}
We summarize the above derivation in the following proposition.
\begin{proposition}
Let $u\in H^1(\Omega)$ be the unique solution of (\ref{eqn:HE}) and 
$y$ the time dependent solution of (\ref{eqn:WE}) corresponding to a minimizer
$(v_0,v_1)\in H^1(\Omega)\times L^2(\Omega)$ of $J$, i.e. $J(v_0,v_1)=0$. 
 Then $u$ is given by (\ref{eqn:shifted-fundamental-frequency-solution}) with $\eta = 0$
if $|\Gamma_S|>0$ or $|\Gamma_D|>0$, and with $\eta$ given by (\ref{eqn:cc_neumann_eta})
 when $\Gamma_N = \partial \Omega$.  
\end{proposition}
Not only does the above filtering approach allow us to use the original cost functional $J$,
it also involves a negligible computational effort or storage amount, 
as the time integral for $\widehat y$ can be calculated cumulatively via numerical 
quadrature during the solution of the wave equation (\ref{eqn:WE}).

\subsection{The CMCG Algorithm \label{algo:CMCG}}
To minimize the quadratic cost functional $J$ defined by (\ref{eqn:cost_functional_J})
over $H^1(\Omega)\times L^2(\Omega)$, 
a natural choice 
is the conjugate gradient (CG) method \cite{BGP1998}, which requires the Fr{\'e}chet derivative of $J$ at $v=(v_0,v_1)$: 
\begin{eqnarray}
  \langle J'(v), \delta v \rangle
  &=&   
  - \int_{\Omega} \nabla (y(x,T)-v_0(x)) \cdot \nabla \delta v_0(x) \ dx  
  \label{eqn:J1_derivative_final}\\
  &&   
  - \int_{\Omega} \frac{1}{c^{2}(x)}(y_t(x,T)-v_1(x))\delta v_1(x) \ dx
  \nonumber\\
  &&   
  +  
   \int_{\Omega}
      \frac{1}{c^{2}(x)}\big(p(x,0)\delta v_1(x)-p_t(x,0) \delta v_0(x))
   \ dx
   \\
&&   
  +
    \intGammaS \hspace*{0.5em}\frac{1}{c(x)}\ p(x,0)\delta v_0(x)
    \ ds   
   \nonumber
.  
\end{eqnarray}
Here
$\delta v= (\delta v_0, \delta v_1)$ denotes an arbitrary perturbation, $\langle\cdot,\cdot\rangle$ the standard duality pairing, and $p$ the solution of the adjoint (backward) wave equation:
\begin{subequations}
  \label{eqn:AWE}  
  \begin{alignat}{4}
    \frac{1}{c^{2}(x)}\frac{\partial^2}{\partial^2 t}p(x,t) \ -\ \laplace p(x,t)
    &\ = \  0 ,&&\qquad x\in\Omega, \ t>0,
    \\
    \dfrac{\partial p(x,t)}{\partial n} -  \dfrac{1}{c}\  \frac{\partial}{\partial t}p(x,t) &\ = \  0 ,&&\qquad x\in\Gamma_S, t>0
    \\
    \dfrac{\partial p(x,t)}{\partial n} &\ = \  0 ,&&\qquad x\in\Gamma_N, t>0,
    \\
    p(x,t) &\ = \  0 ,&&\qquad x\in\Gamma_D, t>0,
    \\
    p(x,T) \ = \ p_0(x),
    \quad
    \frac{\partial p(x,T)}{\partial t} &\ = \ p_1(x)
    ,&&\qquad x\in\Omega,
    \label{eqn:AWE_INIT}
  \end{alignat}
and the initial conditions satisfy for any $w\in H_D^1(\Omega)$
\begin{alignat*}{4}
  p_0(x)
  &\ = \
  y_t(x,T) - v_1(x), 
  \qquad x \in \Omega,
\\   
    \intOmega  
        \frac{p_1(x)}{c^2(x)} w(x)
    \ dx
    &\ = \   
        \intGammaS \hspace*{0.5em} \frac{p_0(x)}{c(x)} w(x)\ ds    
        -
        \intOmega   \nabla (y(x,T)-v_0(x)) \cdot \nabla w(x) \ dx
    .  
\end{alignat*}   
\end{subequations}
The derivation of (\ref{eqn:J1_derivative_final}) and (\ref{eqn:AWE}) can be found in \cite{BGP1998}.
In each CG iteration the derivative $J'(v)$ requires the solution of 
the forward and backward (adjoint) wave equations (\ref{eqn:WE}) and (\ref{eqn:AWE}) 
over one period $[0,T]$. Moreover, each CG iteration requires an explicit (Riesz) representer 
$\tilde g= (\tilde g_0, \tilde g_1)\in H_D^1(\Omega) \times L^2(\Omega)$
of the gradient $g=(g_0,g_1)= J'(v)$ defined in (\ref{eqn:J1_derivative_final}),
which is determined by solving the symmetric and coercive elliptic problem \cite{BGP1998,MS2014}:
\begin{subequations}
  \label{eqn:CMCG_PRECONDITIONIER}
  \begin{alignat}{4}  
    (\nabla \tilde g_0,\nabla \varphi)
    &\ =\ 
    \int_{\Omega}g_0(x)\varphi(x) \ dx      
    \nonumber\\    
    &\ =\ 
    \int_{\Omega} \nabla (v_0(x)-y(x,T)) \cdot \nabla \varphi(x) - \frac{1}{c^2(x)}p_t(x,0)\varphi(x) \ dx
    \nonumber\\    
    &\ \ +     
    \intGammaS \hspace*{0.5em}\frac{1}{c(x)} \ p(x,0)\varphi(x) \ ds
    ,
    \quad \forall \varphi\in H_D^1
    ,    
    \label{eqn:CMCG_PRECOND_STRONGLYELLIPTIC}
\\    
    \tilde g_1
     &\ =\
     g_1
     \ =\    
     v_1-y_t(\cdot,T)
      +
     c^{-2}p(\cdot,0)
     .
  \end{alignat}   
\end{subequations}
For the sake of completeness, we list the full CMCG Algorithm -- see \cite{BGP1998,GT2018}:

\smallskip
\setlist{nolistsep,leftmargin=*}
\noindent{\bf CMCG Algorithm. }
\textit{
 \begin{enumerate}[label=(\arabic*)]
  \item Initialize $v^{(0)}=(v^{(0)}_{0},v^{(0)}_{1})$ (initial guess).
  \item Solve the forward and the backward wave equations (\ref{eqn:WE}) and (\ref{eqn:AWE})
        to determine the gradient of $J$, $g^{(0)} = J'(v^{(0)})$,
        defined by (\ref{eqn:J1_derivative_final}).
  \item Solve the coercive elliptic problem  (\ref{eqn:CMCG_PRECONDITIONIER}) with $g=g^{(0)}$ to determine the 
  new search direction $\tilde g^{(0)}$.
  \item Set $r^{(0)} = d^{(0)} = \tilde g^{(0)}$.
  \item For $\ell=1,2,\ldots$
  \begin{enumerate}[label={\small5.\arabic*\ }]
    \item \label{cmcgalgo:wave} Solve the wave equation (\ref{eqn:WE})     
    with $f=g_D=g_S=g_N=0$ and the initial values $d^{(\ell)}=(d^{(\ell)}_{0},d^{(\ell)}_{1})$ and 
    the backward wave equation (\ref{eqn:AWE}).    
    Compute the gradient $g^{(\ell)} = J'(d^{(\ell)})$ defined by (\ref{eqn:J1_derivative_final}). 
    \item \label{cmcgalgo:elliptic} Solve the coercive elliptic problem  (\ref{eqn:CMCG_PRECONDITIONIER}) with $g=g^{(\ell)}$ to get $\tilde g^{(\ell)}$.
    \item
    $\alpha_{\ell} = 
    \dfrac{
        \|\ \nabla r^{(\ell)}_{0}\|_{L^2(\Omega)}^2
        +
        \|(1/c)\ r^{(\ell)}_{1}\|_{L^2(\Omega)}^2
      }{
          (\ \nabla \tilde g^{(\ell)}_{0},\nabla d^{(\ell)}_{0})_{L^2(\Omega)}
        +
          ((1/c^2)\ \tilde g^{(\ell)}_{1}, d^{(\ell)}_{1})_{L^2(\Omega)}
      }$
    \item $v^{(\ell+1)} = v^{(\ell)}- \alpha_{\ell} d^{(\ell)}$
    \item $r^{(\ell+1)} = r^{(\ell)} - \alpha_{\ell} \tilde g^{(\ell)}$
    \item     
        $\beta_{\ell} = \dfrac{                  
                      \|\nabla r^{(\ell+1)}_{0}\|_{L^2(\Omega)}^2 
                      +
                      \|(1/c)\ r^{(\ell+1)}_{1}\|_{L^2(\Omega)}^2
                    }
                    {                  
                      \|\nabla r^{(\ell)}_{0}\|_{L^2(\Omega)}^2 
                      +
                      \|(1/c)\ r^{(\ell)}_{1}\|_{L^2(\Omega)}^2
                    }
                  $
    \item $d^{(\ell+1)} = r^{(\ell+1)} + \beta_{\ell} d^{(\ell)}$
    \item Stop when the relative residual lies below the given tolerance $tol$
      \begin{equation*}
        \sqrt{
                    \dfrac{
                      \|\nabla r^{(\ell+1)}_{0}\|_{L^2(\Omega)}^2 
                      +
                      \|(1/c)\ r^{(\ell+1)}_{1}\|_{L^2(\Omega)}^2
                    }
                    {                  
                      \|\nabla r^{(0)}_{0}\|_{L^2(\Omega)}^2 
                      +
                      \|(1/c)\ r^{(0)}_{1}\|_{L^2(\Omega)}^2
                    }
                } \leq tol
        .
      \end{equation*}      
  \end{enumerate}
  \item Return approximate solution $u_h$ of (\ref{eqn:HE}) given by
      \[
        \ucmcg  = v_0^{(\ell)} + \frac{i}{\omega} v_1^{(\ell)}.
      \]
 \end{enumerate}
}

Since $\tilde g_0\in H^1(\Omega)$, the updates of $r_0^{(k)}$, $d_0^{(k)}$ and $v_0^{(k)}$
in Steps 5.4, 5.5 and 5.7 in the CMCG Algorithm also remain in $H^1(\Omega)$.
We emphasize that (\ref{eqn:CMCG_PRECOND_STRONGLYELLIPTIC}) is independent of $\omega$
and leads to a symmetric and positive definite linear system, which can be solved efficiently and in parallel
with standard numerical (multigrid, domain decomposition, etc.) methods \cite{DJN2015,BDGST2017}.

\section{Controllability methods for first-order formulations}
The CMCG Algorithm from Section \ref{algo:CMCG} iterates on the initial value 
$(v_0,v_1)\in H^1(\Omega)\times L^2(\Omega)$
of the second-order wave equation (\ref{eqn:WE}) until its solution is $T$-time periodic.
However, the gradient of the cost functional $J(v_0,v_1)$, which is needed during the CG update,
only lies in the dual space $H^{-1}(\Omega)\times L^2(\Omega)$.
To ensure that the solution remains sufficiently regular and in 
$H^1(\Omega)\times L^2(\Omega)$,
the corresponding Riesz representative is computed 
at every CG iteration by solving the strongly elliptic problem (\ref{eqn:CMCG_PRECOND_STRONGLYELLIPTIC}).
In \cite{GR2006}, Glowinski et al. derived 
an equivalent first-order formulation for sound-soft scattering problems,
where the solution instead lies
in $(L^2(\Omega))^{d+1}$,
which is reflexive. 
As a consequence, all CG 
iterates automatically lie in the correct solution space $(L^2(\Omega))^{d+1}$, while
the solution of (\ref{eqn:CMCG_PRECOND_STRONGLYELLIPTIC}) 
is no longer needed.

\subsection{First-order formulation for general boundary conditions}
Again, we always assume for any particular choice of $\omega$, $c(x)$, $f$ and combination of boundary conditions that (\ref{eqn:HE}) has a unique solution $u\in H^1(\Omega)$.
Following \cite{GR2006}, we now let $v=y_t$, ${\bf p}=\nabla y$ and rewrite the time-dependent wave 
equation (\ref{eqn:WE}) in first-order form:
\begin{subequations}  
\label{eqn:mixWE}
\begin{alignat}{4} 
  \frac{1}{c^2(x)}v_{t}(x,t) - \nabla \cdot  {\bf p}(x,t)
    &\ = \ \real{f(x)\e^{-i\omega t}}, &&\quad x\in\Omega, \ t>0,
  \label{eqn:mixWE_INT}
  \\
  \frac{\partial}{\partial t}{\bf p}(x,t) 
    &\ = \ \nabla v(x,t), &&\quad x\in\Omega, \ t>0,
  \label{eqn:mixWE_MIX}
  \\  
  {\bf p}(x,t)\cdot {\bf n} + \frac{1}{c(x)}v(x,t)
    &\ = \  \real{g_S(x)\e^{-i\omega t}}, &&\quad x\in\Gamma_S, t>0,
  \label{eqn:mixWE_ABC}
  \\
  {\bf p}(x,t)\cdot {\bf n}
    &\ = \  \real{g_N(x)\e^{-i\omega t}}, 
    &&\quad x\in\Gamma_N, t>0
,
  \label{eqn:mixWE_NBC}
  \\
  v(x,t) 
    &\ = \ \real{-i\omega g_D(x)\e^{-i\omega t}}, &&\quad x\in\Gamma_D, t>0
  \label{eqn:mixWE_DBC}
\end{alignat}
with the initial conditions
\begin{alignat} {4}
  {\bf p}(x,0) 
    \ = \ {\bf p}_0(x)\in\mathbb{R}^d, \quad
  v(x,0) 
    &\ = \ v_0(x)\in\mathbb{R},&&\qquad x\in\Omega.
\end{alignat}
\end{subequations}
Hence, the solution $({\bf p},v)$ of (\ref{eqn:mixWE}) 
lies in the function space $\mathcal{Q}$~\cite{E2010,P1983},
\begin{equation}
    \label{eqn:mixSpaceTime-FunctionSpace}
    \mathcal{Q} \ =\  C^0([0,T];H(\operatorname{div};\Omega) \times L^2(\Omega)) \cap  C^1([0,T];(L^2(\Omega))^{d+1}).    
\end{equation}
In terms of ${\bf p}$ and $v$,
the energy functional $J$ defined in  (\ref{eqn:cost_functional_J}) now
becomes 
\begin{equation}
  \label{eqn:mix_cost_functional_J}
  \widehat J ({\bf p}_0,v_0) 
  =
    \frac{1}{2} \int_{\Omega}|{\bf p}(x,T) - {\bf p}_0(x)|^2 \ dx
  + \frac{1}{2} \int_{\Omega}\frac{1}{c^{2}(x)}(v(x,T) - v_0(x))^2 \ dx,
\end{equation}
where $({\bf p},v)$ solves (\ref{eqn:mixWE})
with initial value $({\bf p}_0,v_0)\in H(\operatorname{div};\Omega) \times L^2(\Omega)$.

The CMCG Algorithm for the first-order formulation is identical 
to that for the second-order formulation from Section~\ref{algo:CMCG}
except for Steps 2 and 5.1, where $J'$ is now replaced by $\widehat J'$:
\begin{alignat}{4} 
    \label{eqn:mixJ_deriative_2}
    \langle 
      {\widehat J}' ({\bf p}_0,v_0) 
      , 
      (\delta {\bf p}_0, \delta v_0)    
    \rangle  
  =& 
    \int_{\Omega}({\bf p}^*(x,0)-{\bf p}^*(x,T)) \delta {\bf p}_0(x)\ dx
  \\  
  &+  
    \int_{\Omega}(v^*(x,0)-v^*(x,T)) \delta v_0(x)\ dx.
  \nonumber
\end{alignat}
Here $(\delta {\bf p}_0,\delta v_0) \in {\bf P}\times  L^2(\Omega)$ denotes an arbitrary perturbation with
\begin{equation}
    \label{eqn:mixSpaceTime-P}
    {\bf P}=\{{\bf p}\in H(\operatorname{div};\Omega) \ | \  {\bf p}\cdot {\bf n}=0 \text{ on $\Gamma_N$} \},
\end{equation}
whereas $({\bf p}^*,v^*)\in {\bf P}\times L^2(\Omega)$ solves the backward (adjoint)
wave equation in first-order form \cite{GR2006}, that is (\ref{eqn:mixWE})
with $f\equiv g_S\equiv g_N \equiv g_D \equiv 0$ and 
\[
{\bf p}^*(\cdot,T) = {\bf p}(\cdot,T)-{\bf p}_0,
\qquad 
v^*(\cdot,T) = v(\cdot,T)-v_0.
\]

For sound-soft scattering problems ($|\Gamma_D|, |\Gamma_S|>0$),
the functional $\widehat J$ always has a unique (global) minimizer, 
which therefore coincides with the (unique) time-harmonic solution 
$\real{u(x)\e^{-i\omega t}}$ of (\ref{eqn:mixWE}).
For more general boundary value problems, however, the minimizer of $\widehat J$ is not necessarily 
unique, as shown in the following theorem.

\begin{theorem}
\label{thm:mixTP_UNIQUE}
Let $u\in H^1(\Omega)$ be the unique solution of (\ref{eqn:HE}) and 
$({\bf p},v)\in \mathcal{Q}$ be a real-valued solution 
of (\ref{eqn:mixWE}) with initial values $({\bf p}_0,v_0)\in H(\operatorname{div};\Omega) \times L^2(\Omega)$.
If ${\bf p}$ and $v$ are time periodic with period $T=2\pi/\omega$,
then ${\bf p}$ and $v$ admit the Fourier series representation
\begin{subequations}
\label{eqn:mixFourierSeries}
\begin{alignat}{4}
\label{eqn:mixFourierSeries_p}
    {\bf p}(\cdot,t)
    & =  \real{\nabla u \e^{-i\omega t}}        
    + \boldsymbol{\lambda}
    +   \sum_{|\ell|> 1}^{\infty}
                 \boldsymbol{\gamma}_\ell^p\e^{-i\omega \ell t},
\\
\label{eqn:mixFourierSeries_v}
   v(\cdot,t)
    & =
    \omega \imag{u\e^{-i\omega t}}
    +
    \eta
    + 
    \sum_{|\ell|> 1}^{\infty}
        \gamma_\ell^v\e^{-i\omega \ell t},
\end{alignat}
\end{subequations}
where the constant $\eta\in\mathbb{R}$, 
$\boldsymbol{\lambda}\in {\bf P}$ with
\begin{equation}
    \label{eqn:mixElementaryHelmholtzDecomposition} 
    \intOmega \boldsymbol{\lambda} \cdot \nabla \varphi\ dx = 0, \quad \forall \varphi\in H^1(\Omega), \varphi|_{\Gamma_D}\equiv0,
\end{equation}
and the complex-valued eigenfunctions 
$\boldsymbol{\gamma}_\ell^p \in {\bf P}$, $\gamma_\ell^v \in L^2(\Omega)$,  $|\ell|>1$
satisfy
\begin{subequations}  
\label{eqn:mixEigenproblem}
\begin{alignat}{4}
    -c^2(x)\nabla \cdot  {\boldsymbol\gamma}_\ell^p(x)
    +i\omega \ell  {\gamma}_\ell^v(x)
    &\ =\ 0, &&\qquad x\in\Omega, 
    \label{eqn:mixEP_INT}      
  \\
    i\omega\ell  {\boldsymbol\gamma}_\ell^p(x)    
    &\ =\  \nabla {\gamma}_\ell^v(x)
    , 
    &&\qquad x\in\Omega, 
    \label{eqn:mixEP_MIX}
  \\
    c(x){\boldsymbol\gamma}_\ell^p(x) \cdot {\bf n}
    +
    \gamma_\ell^v(x)
    &\ = \  0, &&\qquad x\in\Gamma_S,    
  \label{eqn:mixEP_ABC}
  \\
    {\boldsymbol\gamma}_\ell^p(x) \cdot {\bf n}
    &\ = \  0, &&\qquad x\in\Gamma_N,
  \label{eqn:mixEP_NBC}
  \\
  {\gamma}_\ell^v(x) 
    &\ = \  0, &&\qquad x\in\Gamma_D
  \label{eqn:mixEP_DBC}  
  .
\end{alignat}
\end{subequations}
Furthermore, if $|\Gamma_S\cup \Gamma_D|>0$, then $\eta = 0$.
\end{theorem}
\proof Let 
\begin{equation*}
    {\bf q}(\cdot,t) =  {\bf p}(\cdot,t) - \real{\nabla u \e^{-i\omega t}},
    \quad
    w(\cdot,t)  =  v(\cdot,t) - \omega\imag{u\e^{-i\omega t}}.
\end{equation*}
Then $w$ and ${\bf q}$ satisfy (\ref{eqn:mixWE})
with $f\equiv g_D\equiv g_S\equiv g_N\equiv 0$ and initial values 
\begin{equation*}
{\bf q}(x,0) = {\bf p}_0(x) - \real{\nabla u(x)}
,\quad
w(x,0) = v_0(x) - \omega\imag{u(x)}
,
\quad x\in\Omega.
\end{equation*}
Since ${\bf p}$ and $v$ are $T$-periodic, so are
${\bf q}$ and $w$.
Moreover, the mappings 
\begin{equation*}
t \mapsto ({\bf q}(\cdot,t),\boldsymbol{\psi})
,\qquad 
t \mapsto (w(\cdot,t),\varphi)
\end{equation*}
are $T$-periodic and continuous
for any  $(\boldsymbol{\psi},\varphi) \in {\bf P}\times L^2(\Omega)$~\cite{E2010}.
Hence, they admit the Fourier series representation,
\begin{equation*}
    ({\bf q}(\cdot,t),\boldsymbol\psi)
    = 
    \hspace*{-.25em}\sum_{\ell=-\infty}^{\infty}\hspace*{-.25em}\widehat{\boldsymbol\gamma}_\ell^p\e^{i \omega\ell t}
    ,
    \qquad
    (w(\cdot,t),\varphi)
    =  
    \hspace*{-.25em}\sum_{\ell=-\infty}^{\infty}\hspace*{-.25em}{\widehat \gamma}_\ell^v\e^{i \omega\ell t}    
    ,
\end{equation*}
where
$\boldsymbol{\gamma}_\ell^p\in\mathbb{C}^{d}$
, 
$\widehat{\gamma}_\ell^v \in\mathbb{C}$.
Next, we define
\begin{equation}  
\label{eqn:FourierSeries_Coefficients}
\boldsymbol{\gamma}_\ell^p(x)= \frac{1}{T}\intT  {\bf q}(x,t)\e^{-i \omega\ell t}dt, 
\qquad
{\gamma}_\ell^v(x)=  \frac{1}{T}\intT  w(x,t)\e^{-i \omega\ell t} \ dt, 
\end{equation}
which implies that 
\begin{equation*} 
    \widehat{\boldsymbol\gamma}_\ell^p \ =\  (\boldsymbol{\gamma}_\ell^p,\boldsymbol{\psi})
    ,
    \qquad
    \widehat{\gamma}_\ell^v \ =\  ({\gamma}_\ell^v,\varphi)
    .
\end{equation*}

We shall now show that $\boldsymbol{\gamma}_\ell^p$ and $\gamma_\ell^v$
satisfy (\ref{eqn:mixEigenproblem}) for all $|\ell|\ge1$.
First, integration by parts, 
(\ref{eqn:mixWE_INT})-(\ref{eqn:mixWE_MIX}) and the periodicity of ${\bf q}$ and $w$ imply 
\begin{alignat*}{4}
    {\gamma}_\ell^v (x)
     & =        
    \frac{1}{T}\intT w_t(x,t)\frac{\e^{-i\omega\ell t}}{i\omega \ell } \ dt    
     -
    \frac{w(x,t)\e^{-i\omega\ell t}}{i\omega \ell  T} \bigg|_{0}^{T}
     =     
    \frac{1}{T}\intT c^2(x)\nabla \cdot {\bf q}(x,t)
    \frac{\e^{-i\omega\ell t}}{i\omega\ell} \ dt
    ,
    \\
    \boldsymbol{\gamma}_\ell^p (x)
     & =         
    \frac{1}{T}\intT {\bf q}_t(x,t)\frac{\e^{-i\omega\ell t}}{i\omega \ell } \ dt    
    -
    \frac{{\bf q}(x,t)\e^{-i\omega\ell t}}{i\omega \ell  T} \bigg|_{0}^{T}
     =     
    \frac{1}{T}\intT \nabla w(x,t)
    \frac{\e^{-i\omega\ell t}}{i\omega\ell} \ dt
    .
\end{alignat*}
Together with definition (\ref{eqn:FourierSeries_Coefficients}) of ${\boldsymbol\gamma}_\ell^p$ and $\gamma_\ell^{v}$, we thus immediately obtain
\begin{equation*}
    i\omega \ell{\gamma}_\ell^v - c^2\nabla  \cdot {\boldsymbol\gamma}_\ell^p
    \ =\ 0    
    ,
    \quad     
    i\omega\ell\boldsymbol{\gamma}_\ell^p
    \ =\ 
    \nabla {\gamma}_\ell^v
    \qquad    
    \text{in $\Omega$}
    .
\end{equation*}
Since $w(x,t)=0$ for $x\in\Gamma_D$, we infer from (\ref{eqn:FourierSeries_Coefficients}) that 
\[
    \intGammaD \gamma_\ell^v(x) \varphi (x) \ ds
    =  
    \frac{1}{T}\intT \intGammaD w(x,t) \varphi(x)\ ds \e^{-i\omega\ell t} \ dt  
     = 
    0    
    ,
    \qquad
    \varphi \in L^2(\Gamma_D),
\]
and hence $\gamma_\ell^v$ satisfies (\ref{eqn:mixEP_DBC}).
Similarly, (\ref{eqn:mixEP_ABC}), (\ref{eqn:mixEP_NBC}) follow from 
the fact that  ${\bf q}$ and $w$  satisfy
(\ref{eqn:mixWE_ABC}), (\ref{eqn:mixWE_NBC})
with $g_N\equiv g_S \equiv 0$.
Hence $\boldsymbol{\gamma}_\ell^p$, $\gamma_\ell^v$ satisfy (\ref{eqn:mixEigenproblem}) 
for all $|\ell|\ge1$.
In fact for $\ell=1$, (\ref{eqn:mixEigenproblem}) corresponds
to (\ref{eqn:HE}) in first-order formulation with $\boldsymbol{\gamma}_{1}^{p}=\nabla \overline{u}$, $\gamma_{1}^{v}=i\omega \overline{u}$,
homogeneous boundary conditions and no sources. By uniqueness,
$\boldsymbol{\gamma}_{1}^{p}$ and $\gamma_{1}^v$, together with their complex conjugates, 
are therefore identically zero.

Next, we consider $\boldsymbol{\gamma}_0^{p}$, $\gamma_0^{v}$.
Again, since ${\bf q}$ and $w$ satisfy (\ref{eqn:mixWE_INT})-(\ref{eqn:mixWE_DBC})
with $f=0$ and homogeneous boundary conditions, we obtain from 
(\ref{eqn:FourierSeries_Coefficients}) with $\ell=0$ and the periodicity of ${\bf q}$ and $w$
\begin{alignat}{4}
  \label{eqn:THM_MIXED_LAMBDA_1}
 \intOmega (\nabla \cdot \boldsymbol{\gamma}_0^p) \ \varphi \ dx
 &=
  \frac{1}{T}\intT\intOmega \frac{1}{c^{2}}{w}_t \varphi  \ dx dt
  = 
  0,    
  \qquad \forall \varphi\in L^2(\Omega),
\\
  \label{eqn:THM_MIXED_LAMBDA_2}
 \intOmega {\gamma}_0^v  \nabla \cdot \boldsymbol{\psi} \ dx
  &=
\frac{1}{T}\intT\intOmega {\bf q}_t \cdot \boldsymbol{\psi} \ dx dt
  -
  \frac{1}{T}\intT\intGammaS w \ \boldsymbol{\psi} \cdot {\bf n}\ ds dt
  \\
  &=  
  \frac{1}{T}\intT\intGammaS c\ {\bf q} \cdot {\bf n} \ \boldsymbol{\psi}\cdot {\bf n}\ ds dt
  =
  \intGammaS c\ \boldsymbol{\gamma}_0^{p} \cdot {\bf n}\ \boldsymbol{\psi}\cdot {\bf n}\ ds,
  \qquad \forall \boldsymbol{\psi}\in {\bf P}
  \nonumber
  . 
\end{alignat}
In particular, (\ref{eqn:THM_MIXED_LAMBDA_1})-(\ref{eqn:THM_MIXED_LAMBDA_2}) 
implies with $\varphi=\gamma_0^{v}$ and $\boldsymbol{\psi}=\boldsymbol{\gamma}_0^{p}$
that 
\begin{equation*}
    \intGammaS c |\boldsymbol{\gamma}_0^{p} \cdot {\bf n}|^2\ ds
    = 0,  
\end{equation*}
and hence, $\boldsymbol{\gamma}_0^{p} \cdot {\bf n}=0$ on $\Gamma_S$, since $c>0$.
Moreover, Green's formula, together with (\ref{eqn:THM_MIXED_LAMBDA_1})  and the homogeneous boundary conditions, implies that 
\begin{alignat*}{4}
  \intOmega \boldsymbol{\gamma}_0^p \cdot  \nabla \varphi \ dx
  &=
  -
  \intOmega (\nabla\cdot \boldsymbol{\gamma}_0^p) \  \varphi \ dx
  +
  \intPartialOmega \boldsymbol{\gamma}_0^p \cdot {\bf n} \ \varphi \ ds 
  =
  0,
  \quad \forall\varphi \in H^1(\Omega), \varphi|_{\Gamma_D}\equiv 0,
\end{alignat*}
and therefore $\boldsymbol{\lambda}=\boldsymbol{\gamma}_0^{p}$ satisfies (\ref{eqn:mixElementaryHelmholtzDecomposition}).

To show that $\gamma_0^v$ is constant, we now let $\varphi\in \mathcal{C}_c^\infty(\Omega)$ and 
$\boldsymbol{\psi} = {\bf e}_j\varphi \in H(\operatorname{div};\Omega)$, $j=1,\ldots,d$, where ${\bf e}_j$ is  
the $j$-th unit basis vector of $\mathbb{R}^d$.
Integration of (\ref{eqn:mixWE_MIX}) over $[0,T]$, definition (\ref{eqn:FourierSeries_Coefficients}) with $\ell=0$ and the periodicity of  ${\bf q}$ then
yield 
\begin{alignat}{4}
  \label{eqn:THM_MIXED_ETA_1}
   0 
   & =  \frac{1}{T}\intT\intOmega {\bf q}_t \ \boldsymbol{\psi}\ dx dt
   \ =\ - \frac{1}{T}\intT\intOmega w \ \nabla \cdot \boldsymbol{\psi}\ dx dt 
   \ =\ -  \intOmega \gamma_0^v \ \frac{\partial \varphi}{\partial x_j}\ dx.
\end{alignat}
From (\ref{eqn:THM_MIXED_ETA_1}), we conclude that $\partial_{x_j} \gamma_0^v = 0$, $j=1,\ldots,d$,
which implies
\begin{equation*}
  \gamma_0^v(x) \equiv \eta, \qquad \gamma_0^v\in H^1(\Omega).
\end{equation*}

Since $\gamma_0^v$ satisfies (\ref{eqn:mixWE_DBC}) with $\ell=0$, $\eta=\gamma_0^v = 0$, if $|\Gamma_D|>0$.
Similarly, if $|\Gamma_S>0|$,  (\ref{eqn:mixWE_ABC}), together with $\boldsymbol{\gamma}_0^{p} \cdot {\bf n}=0$ on $\Gamma_S$, yields
\begin{equation*}
  0 
  =
  \frac{1}{T}\intT \big[c(x) {\bf q}(x,t) \cdot {\bf n}  + w(x,t)  \big]\ dt
  =
  c(x)\boldsymbol{\gamma}_0^p(x) \cdot {\bf n}  + {\gamma}_0^v(x)   
  =
  \eta,  
  \;\;
  x\in\Gamma_S.
\end{equation*}
Thus, $\eta = 0$ when $|\Gamma_D\cup \Gamma_S|>0$, which completes the proof.

\qed

For sound-soft scattering problems, where $|\Gamma_D|>0$ and $|\Gamma_S|>0$, $\eta=0$ and all eigenfunctions $\boldsymbol{\gamma}_\ell^p$, $\gamma_\ell^v$, $|\ell|>1$
of (\ref{eqn:mixEigenproblem}) trivially vanish in (\ref{eqn:mixFourierSeries})
\cite{CX2006}. Therefore, (\ref{eqn:mixFourierSeries})-(\ref{eqn:mixElementaryHelmholtzDecomposition}) 
in Theorem \ref{thm:mixTP_UNIQUE} with $t=0$ 
imply that 
\begin{eqnarray*}
({\bf p}_0,\nabla \varphi)
&=&
     (\real{\nabla u},\nabla \varphi)
     ,
     \qquad \varphi\in H^1(\Omega), \; \varphi|_{\Gamma_D}=0
,
\\     
v_0
&=&
     \omega\imag{u}
.
\end{eqnarray*}
From the real part of (\ref{eqn:HE}) we than conclude that	
\begin{equation}
    \label{eqn:mix_solution}
    u = -k^{-2}\big(\real{f} + \nabla \cdot {\bf p}_0) + i\omega^{-1}v_0.
\end{equation}

\subsection{Fundamental frequency filtering for first-order formulation\label{sec:mixFilterFundamentalFrequency}}
When the CMCG method is applied to the first-order formulation (\ref{eqn:mixWE}),
any minimizer of $\widehat J({\bf p}_0,v_0)=0$ generally consists of 
spurious perturbations $\eta$, $\boldsymbol{\lambda}$ and 
eigenfunctions $\boldsymbol{\gamma}_\ell^p$, $\gamma_\ell^v$ superimposed on the 
desired (unique) solution $u$ of (\ref{eqn:HE}).
To extract $u$ from $({\bf p}_0,v_0)$, we apply a filtering approach, 
similar to that in Section \ref{sec:FilterFundamentalFrequency}, and thereby
restore uniqueness. 
Again, we multiply the Fourier series representation in (\ref{eqn:mixFourierSeries})
of $v$ by $\e^{i\omega t}$ and integrate over $(0,T)$.
Since $\eta$ and $\boldsymbol{\lambda}$ are independent of time, while $\e^{i\omega t}$ is orthogonal
to $\e^{i\omega \ell t}$, $|\ell|>1$, all spurious modes vanish
and the resulting expression simplifies to:
\begin{alignat*}{4}
    \frac{2}{T}\int_{0}^{T}v(\cdot,t)\e^{i\omega  t}\ dt   
    &\ =\          
    \frac{2}{T}\int_{0}^{T}\real{-i\omega u\e^{-i\omega  t}}\e^{i\omega  t}\ dt
    &&\ =\          
    -i\omega u
    ,
\end{alignat*}
which immediately yields
\begin{equation}
    \label{eqn:TH_FILTERED_MIXED_SOLUTION}
    u(x)
    = 
    \frac{2i}{T\omega}\int_{0}^{T}v(\cdot,t)\e^{i\omega  t}\ dt
    .
\end{equation}
We summarize this result in the following proposition.
\begin{proposition}
 Let $u\in H^1(\Omega)$ be the unique solution of (\ref{eqn:HE})
 and $({\bf p},v)\in \mathcal{Q}$ be a $T$-time periodic solution of (\ref{eqn:mixWE}).
 Then $u$ is  given by (\ref{eqn:TH_FILTERED_MIXED_SOLUTION}) .
\end{proposition}

\subsection{Hybrid DG FE-Discretization}
In \cite{GR2006}, Glowinski et al. used standard Raviart-Thomas (RT) finite elements 
to discretize (\ref{eqn:mixWE}).
Since no mass-lumping is available  
for RT elements
on triangles or tetrahedra~\cite{BJT2000}, 
each time-step then
requires the inversion of the mass-matrix.
To avoid that extra computational cost, which strongly impedes parallelization,
we instead consider the recent hybrid discontinuous Galerkin (HDG) FEM \cite{CNPS2016}
to discretize (\ref{eqn:HE}) in its corresponding first-order formulation together with (\ref{eqn:mixWE}).
Then, the mass-matrix is block-diagonal, 
with (small and constant) block size equal to the number of 
dof's per element, so that the 
time-stepping scheme becomes truly explicit and inherently
parallel.

Let $\mathcal{T}_h$ denote a regular triangulation of $\Omega_h$,
$\mathcal{E}_h$ the set of all faces and $\mathcal{P}^r$ the space of polynomials of degree $r$. In addition, we define 
\begin{alignat}{4}
{\bf P}_h \ =\ & \{{\bf r}\in (L^2(\Omega))^d: {\bf r}|_K \in (\mathcal{P}^r(K))^d, \forall K \in \mathcal{T}_h \}
,
\\
V_h \ =\ & \{w\in L^2(\Omega): w|_K \in \mathcal{P}^r(K), \forall K \in \mathcal{T}_h \}
,
\\
M_h \ =\ & \{\mu \in L^2(\mathcal{E}_h) : \mu|_F \in \mathcal{P}^r(F), \forall F \in \mathcal{E}_h \}
.
\end{alignat}

\begin{subequations}
Following \cite{CNPS2016},
the HDG Galerkin FE formulation reads: \\ Find  $({\bf p}_h,v_h,\widehat v_h)\in {\bf P}_h \times V_h \times M_h$  
such that 
\label{eqn:HDG_wave_equation}
\begin{alignat}{4}    
    \big( \frac{\partial {\bf p}_h}{\partial t}, {\bf r}\big)_{K}
    \ =\ &
    -
    (v_h, \nabla \cdot {\bf r})_{K}
    +
    \langle\widehat v_h, {\bf r}\cdot {\bf n} \rangle_{\partial K}  
    ,
    \\
    \big(\frac{1}{c}\ \frac{\partial v_h}{\partial t}, w \big)_{K}
    \ =\ &
    (f, w)_{K}
    -
    ({\bf p}_h, \nabla w)_{K}
    +
    \langle \widehat {\bf p}_h\cdot {\bf n}, w \rangle_{\partial K}
    ,
    \\    
    \langle    
    \widehat {\bf p}_h\cdot {\bf n}
    +
    \frac{1}{c}\widehat v_h
    ,    
    \mu
    \rangle_{\partial K \cap \Gamma_S}  
    \ =\ &
    \langle    
    \real{g_S \exp(-i \omega t)},
    \mu
    \rangle_{\partial K \cap \Gamma_S}
    ,
    \\
    \langle    
    \widehat {\bf p}_h\cdot {\bf n}    
    ,    
    \mu
    \rangle_{\partial K \cap \Gamma_N}  
    \ =\ &
    \langle    
    \real{g_N \exp(-i \omega t)},
    \mu
    \rangle_{\partial K \cap \Gamma_N}  
    ,
    \\
    \langle    
    \widehat v_h
    ,
    \mu
    \rangle_{\partial K \cap \Gamma_D}  
    \ =\ &
    \langle    
    \omega\imag{g_D \exp(-i \omega t)},
    \mu
    \rangle_{\partial K \cap \Gamma_D}
    ,
 \end{alignat} 
for all $({\bf r},w,\mu)\in {\bf P}_h \times V_h  \times M_h$, $K \in \mathcal{T}_h$ and $t\in[0,T]$,
where $(\cdot, \cdot)_K$ and $\langle \cdot, \cdot \rangle_D$ denote the $L^2$-inner product on $K$ or $D$, respectively and the numerical flux is
\begin{alignat}{4}
    \widehat {\bf p}_h \cdot {\bf n}  
    = 
    {\bf p}_h \cdot {\bf n}
    -
    \tau (v_h-\widehat{v}_h)
    \quad \text{on $\partial K$}
\end{alignat} 
with the stabilization function $\tau$ from \cite{CNPS2016}.
In addition, $({\bf p_h}, v_h)$ satisfies the initial conditions
\begin{alignat}{4} 
{\bf p}_h (x,0) \ =\ & {\bf p}_0(x)
,
\qquad
v_h (x,0) \ =\ & v_0(x)
,
\qquad
x\in\Omega.
\end{alignat} 
\end{subequations}
For the time integration of (\ref{eqn:HDG_wave_equation}), we use the
standard explicit fourth-order Runge-Kutta (RK4) method.

\subsection{Convergence and superconvergence}
\begin{figure}[t] 
 \centering
    \includegraphics{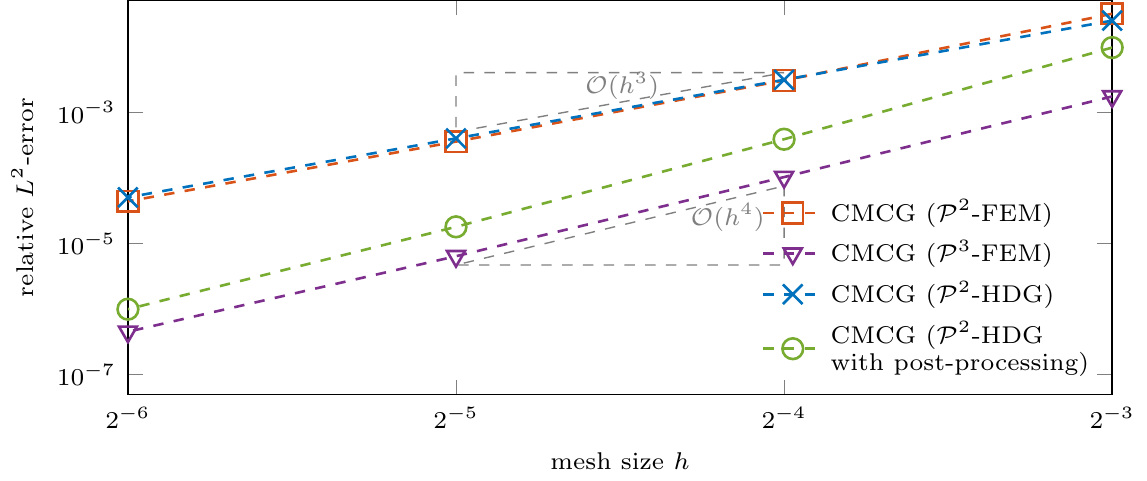}
    \caption{Convergence and superconvergence:
    the numerical error $\|u-u_h\|$ vs. mesh size $h=2^{-i}$, $i=3,\ldots,6$,
    obtained 
    with the CMCG method for the second-order formulation with $\mathcal{P}^{2}$-/$\mathcal{P}^{3}$-FE
    or
    for the first-order formulation with $\mathcal{P}^{2}$-HDG discretization,
    either with or without post-processing.
    }
    \label{fig:hdg-convergence}
\end{figure}
For a FE discretization with piecewise polynomials of degree $r$,
we usually expect convergence as $\mathcal{O}(h^{r+1})$
with respect to the $L^2$-norm.
For the above HDG discretization, however,
an extra power in $h$ can be achieved by applying a cheap
local post-processing step~\cite{CNPS2016}.
The same  (super-) convergence in space of 
order $r+2$ using only $P^r$-FE can be achieved with the CMCG method 
by applying the post-processing step to the numerical solutions $({\bf p}_h^{n_T}, v_h^{n_T})$
of (\ref{eqn:mixWE}) at the final time $T=n_T\Delta t$.

\begin{subequations}
Let $({\bf p}_h^m,v_h^m,{\hat v}_h^m)$ denotes the fully discrete solution
of (\ref{eqn:HDG_wave_equation}) at $t_m=m\Delta t$.
First, we compute the new (more accurate) approximation ${\bf p}_h^{n_T,*}$ of ${\bf p}(\cdot,T)$ by solving the local problem
\begin{alignat*}{4}
    ({\bf p}_h^{n_T,*}, \boldsymbol{\psi})_{L^2(K)}
    &\ = \ 
    -(v_h^{n_T}, \nabla \cdot \boldsymbol{\psi})_{L^2(K)}
    +\langle {\widehat v}_h^{n_T}, \boldsymbol{\psi} \cdot {\bf n} \rangle_{\partial K}
    ,
    \quad 
    \forall \boldsymbol{\psi}\in {\bf P}_h
\end{alignat*}
on each $K\in\mathcal{T}_h$.
Then, we calculate the additional approximations
$y_h^{n_T,*}$ of $y(\cdot, T)=\operatorname{Re}\{u(x)\}$ given by (\ref{eqn:mix_solution}),
$v_h^{n_T,*}$ of $v(\cdot, T)$
in $\mathcal{P}^{r+1}(K)$, which satisfy
\begin{alignat*}{4}
    (\nabla y_h^{n_T,*},\nabla \varphi)_{L^2(K)}
    &\ = \  ({\bf p}_h^{n_T}, \nabla \varphi)_{L^2(K)}, 
    \qquad 
    &\forall \varphi\in\mathcal{P}^{r+1}(K),
    \\
    (y_h^{n_T,*},1)_{L^2(K)} 
    &\ = \ (y_h^{n_T},1)_{L^2(K)},
    \\
    (\nabla v_h^{n_T,*},\nabla \varphi)_{L^2(K)} 
    &\ = \   ({\bf p}_h^{n_T,*}, \nabla \varphi)_{L^2(K)}, 
    \qquad
    &\forall \varphi\in\mathcal{P}^{r+1}(K),
    \\
    (v_h^{n_T,*},1)_{L^2(K)} 
    &\ = \   (v_h^{n_T}, 1)_{L^2(K)},    
\end{alignat*}
\end{subequations}
for any element $K \in \mathcal{T}_h$. 
The new approximate solution  $u$  is then
given by (\ref{eqn:TH_FILTERED_MIXED_SOLUTION})
with ${\bf p}$  and $v$ replaced by
${\bf p}_h^{n_T,*}$ and $v_h^{n_T,*}$.

To illustrate the accuracy and verify the expected convergence rates
for the various FE discretizations in the CMCG method,
we now consider the following one-dimensional solution 
\[
    u(x) = -\exp(ikx)
\]
of (\ref{eqn:HE}) in $\Omega = (0,1)$ 
with $c=1$, $k=5\pi/4$, 
$\Gamma_D=\{0\}$ and $\Gamma_S=\{1\}$.
Figure \ref{fig:hdg-convergence}
shows the error $\|u-u_h\|$
obtained with the CMCG method
for the first-order formulation (\ref{eqn:WE}) and a $P^2$-HDG 
discretization on a sequence of increasingly finer meshes 
$h=2^{-i}$, $i=3,\ldots,6$. Clearly as we refine the mesh, we always
reduce the time-step in the RK4 method to satisfy the CFL stability condition.
The CG iteration stops once the tolerance $tol=10^{-12}$ is reached.
We also compare the solutions obtained with the CMCG method 
applied to the second-order formulation using
a  (continuous) $\mathcal{P}^{2}$ or $\mathcal{P}^{3}$-FEM.
All numerical solutions display the expected optimal convergence of order $r+1$ with polynomials of degree $r$, while the first-order HDG approach even achieves
superconvergence of order $r+2$,
once local post-processing is applied to the final CG iterate.

\subsection{Physically bounded domain \label{sec:mix_hepbd}}
\begin{figure}[t] 
\centering
  \includegraphics{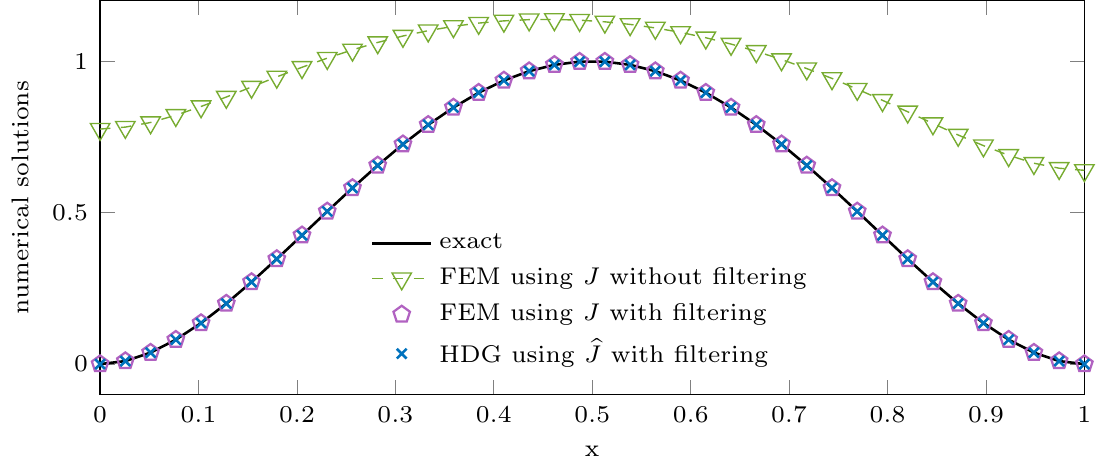}  
  \caption{Physically bounded domain:
  comparison of the exact solution $u$ of (\ref{eqn:HE})
  with the numerical solutions $u_h$ obtained with
  the CMCG method either applied to
  the second-order formulation
  with standard FEM or to the 
  first-order formulation with an HDG discretization.
  }
  \label{fig:one-dimensional_helmholtzneumann}
\end{figure}
In the absence of Dirichlet or impedance boundary conditions, 
the first-order formulation does not yield the 
correct minimizer of $J$. As a simple remedy, we proposed in Section \ref{sec:mixFilterFundamentalFrequency} a filtering procedure which removes the unwanted spurious modes. To illustrate the effectiveness of the filtering procedure, 
we now consider
the exact solution of (\ref{eqn:HE})
\begin{equation}
  u(x) = 16x^2(x-1)^2
\end{equation}
in $\Omega=(0,1)$ with homogeneous Neumann boundary conditions 
and $k=\omega=\pi/4$, $c=1$.
Note that 
$k^2$ is not an eigenvalue of (\ref{eqn:Eigenproblem})
and therefore the solution of (\ref{eqn:HE}) is well-posed.
However, as $(4k)^2=\pi^2$ indeed corresponds to the first eigenvalue 
of the negative Laplacian,
the CMCG method in general
will not yield the correct (unique) solution -- see Theorems \ref{thm:TP_UNIQUE} and \ref{thm:mixTP_UNIQUE}.
Indeed as shown in Figure \ref{fig:one-dimensional_helmholtzneumann},
the original CMCG method \cite{BGP1998}
applied to the second-order formulation with the energy functional $J$ in (\ref{eqn:cost_functional_J})
does not yield the exact solution of (\ref{eqn:HE}),
unlike the numerical solutions obtained after filtering -- see Sections \ref{sec:FilterFundamentalFrequency} and \ref{sec:mixFilterFundamentalFrequency}.

\section{Numerical results}
Here we present a series of numerical examples that illustrate the accuracy, convergence behavior and parallel
performance of the CMCG method.
First, we verify that the numerical solution $\ucmcg$
of (\ref{eqn:HE}) obtained with the CMCG method converges to the numerical solution $\uddm$ 
obtained with a direct solver for the same spatial FE discretization as the time step $\Delta t\rightarrow 0$
in the numerical integration of (\ref{eqn:WE}).
Next, we 
evaluate different stopping criteria for 
the CG iteration in the CMCG Algorithm  from Section~\ref{algo:CMCG}.
We also compare the CMCG Algorithm to a long-time solution of the wave equation without controllability (``do-nothing'' approach)
to demonstrate its effectiveness, in particular for nonconvex obstacles. 
Moreover, we show how an initial run-up
yields a judicious initial guess $(v_0,v_1)$
for the CG iteration thereby further accelerating convergence.
Finally, we apply the CMCG method to large scale scattering problems on a massively parallel architecture,
where the elliptic problem (\ref{eqn:CMCG_PRECONDITIONIER})
is solved in parallel with domain decomposition methods.

\subsection{Semi-discrete convergence \label{sec:semi-discrete-solution}}
\begin{figure}[t]
\centering
  \includegraphics{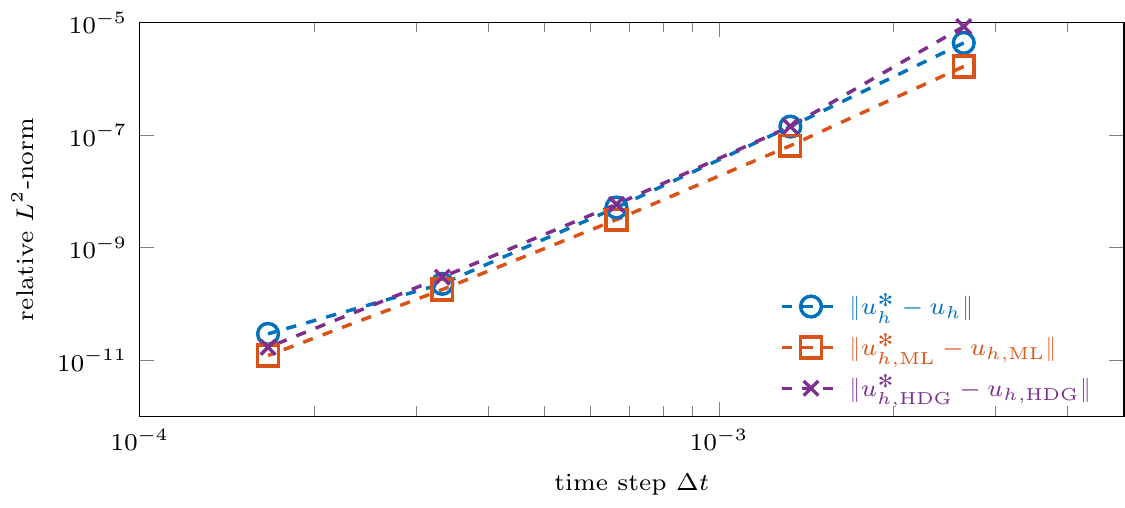}
\caption{Semi-discrete convergence:
Comparison of the numerical solution $\ucmcg$,
obtained with the CMCG method,
and $\uddm$, obtained with a direct solver 
for the same fixed $P^2$-FE discretization ($H^1$-conforming or HDG), both
either with or without mass-lumping (ML)
}
\label{fig:discretization_in_time}
\end{figure}
First, we consider a simple 1D example
to show for a fixed FE-mesh 
that the numerical solution $\ucmcg$, obtained with the CMCG method,
converges to the numerical solution $\uddm$, obtained with a direct solver,
as $\Delta t\rightarrow 0$.
Hence we consider the following solution $u$ of (\ref{eqn:HE})  in $\Omega=(0,1)$
with $\omega=k=6\pi$, $c=1$ and $f\equiv 0$:
\[
    u(x) = \exp(ikx),
    \qquad
    \text{with}
    \qquad
    u(0)=1,
    \qquad
    u'(1)-ik\ u(1) = 0.
\]
Now, let $\uddm(x)$ be the FE Galerkin solution corresponding to the direct solution of the linear system
\begin{equation}
   {\bf A}_h{\bf u}_h^* = {\bf b}_h,
   \label{eqn:discrete_HE}
\end{equation}
resulting from the same
standard $H^1$-conforming or HDG $P^2$-FE discretization of the Helmholtz equation (\ref{eqn:HE})
in second- or first-order formulation, respectively.
For the time integration of (\ref{eqn:WE}) or (\ref{eqn:mixWE})
in the CMCG Algorithm, we use the standard explicit fourth order Runge-Kutta (RK4) method.

Usually we avoid inverting the mass-matrix at each time step via
order preserving mass-lumping \cite{CJRT2001} which, however, introduces
an additional spatial discretization error.
Here to ensure a consistent comparison, we thus compute
$\ucmcg$
and
$\uddm$
both either with, or without, mass-lumping (ML).
For the CG iteration, we always choose $v_0^{(0)}\equiv 0$, $v_1^{(0)}\equiv 0$ 
and fix the tolerance to $tol=10^{-14}$ 
to ensure convergence to machine precision accuracy.

In Figure \ref{fig:discretization_in_time}, we monitor 
the difference between the numerical solution
$\uddm$ or $\uddmHDG$ of (\ref{eqn:discrete_HE}),
obtained with a direct solver, and $\ucmcg$ or $\ucmcgHDG$, obtained 
with the CMCG method using either the second or the first order formulation, respectively.
As expected, for increasingly smaller $\Delta t$ and a fixed 
stringent tolerance in the CG iteration, 
the numerical solution of the CMCG method always
converges to the discrete solution of the Helmholtz equation
for the same FE discretization.

\subsection{CG iteration and initial run-up \label{sec:square-shaped-scattering-problem}}
\begin{figure}[t]
    \centering
    \captionsetup[subfigure]{justification=centering}
    \hspace*{-0.25cm}     
    \subfloat[convex obstacle \protect\\($145'924$ $P^2$-FE) \label{fig:square-shaped-obstacle-mesh-a}]{%
        \includegraphics[width=4cm]{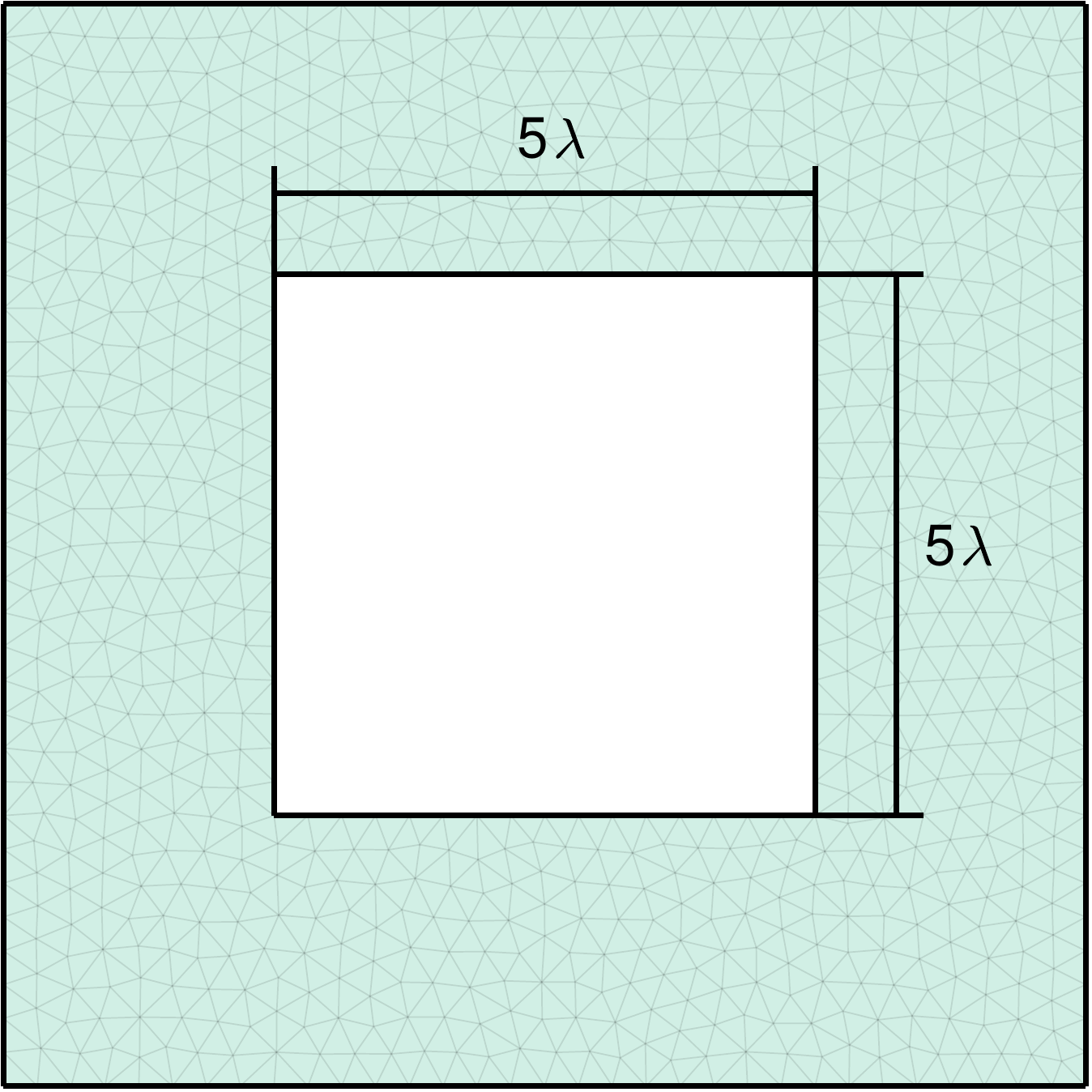}}
    \hspace*{1.5cm}
    \subfloat[nonconvex obstacle \protect\\($176'018$ $P^2$-FE) \label{fig:square-shaped-obstacle-mesh-b}]{%
        \includegraphics[width=4cm]{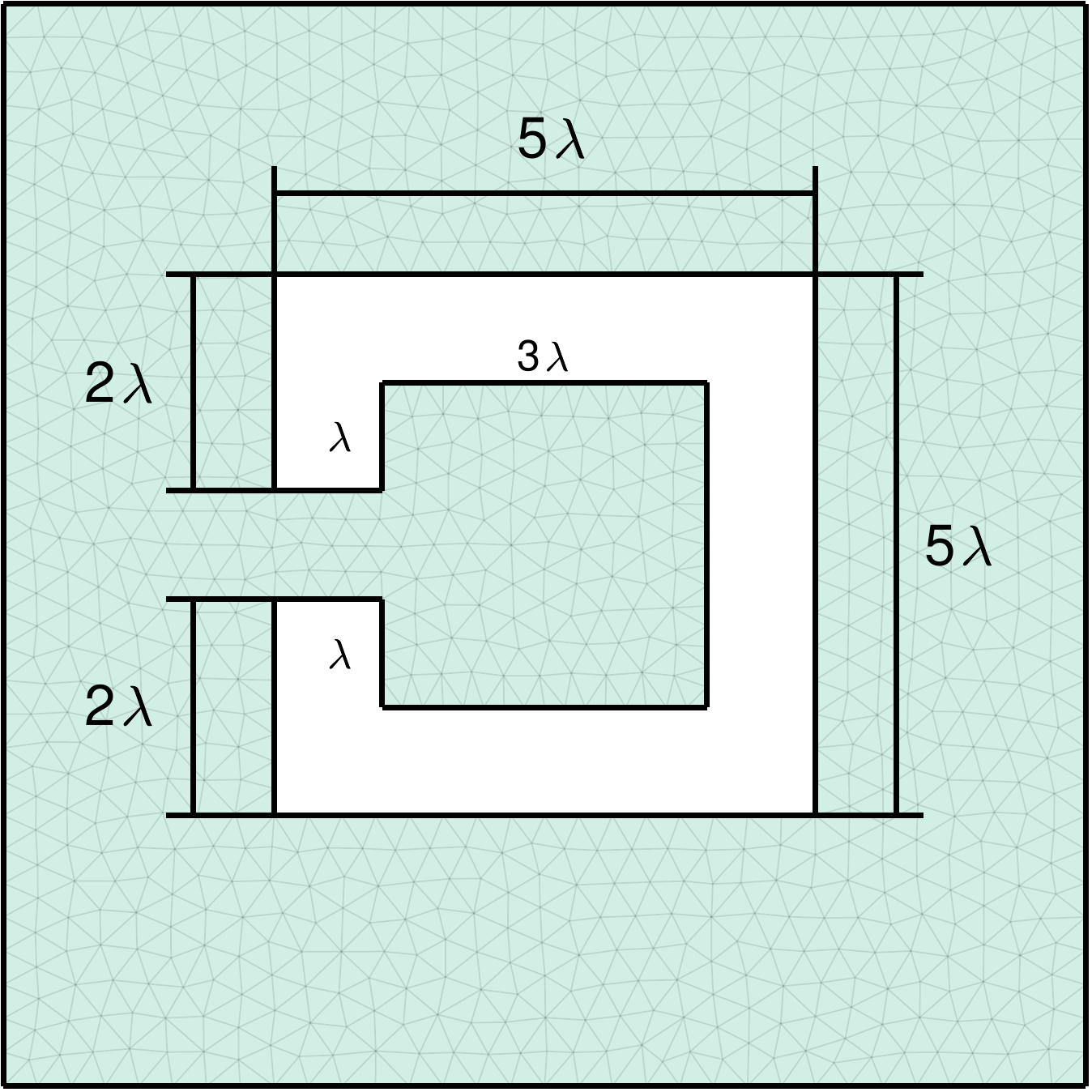}}
    \caption{Computational domain $\Omega$ with a convex square (a) or a nonconvex cavity (b) shaped obstacle}
    \label{fig:square-shaped-obstacle-mesh}
\end{figure}
Next, we first compare different stopping criteria for 
the CG iteration in the CMCG Algorithm
applied to the original second-order 
formulation from Section \ref{sec:2nd-order-formulation}. 
We then illustrate
how the CMCG method greatly accelerates 
the convergence of a solution of the wave equation
to its long-time asymptotic limit, in particular for nonconvex obstacles.
Finally, we show how an initial run-up yields a judicious initial guess 
for the CG iteration, which further accelerates the convergence of the CMCG Algorithm.

Hence, we consider a two-dimensional 
sound-soft scattering problem (\ref{eqn:HE})
with $c\equiv 1$,  $k =\omega=2\pi$, $f\equiv g_D\equiv g_N\equiv 0$ and $g_S=-(\partial_n - ik)u^{in}$
in a bounded square domain $\Omega=(0,10\lambda) \times (0,10\lambda)$, $\lambda=1$,
either with a convex obstacle 
or a  semi-open square shaped cavity.
On the boundary $\Gamma_D$ of the obstacle, we impose a homogeneous Dirichlet condition 
and on the exterior boundary $\Gamma_S$
a Sommerfeld-like absorbing condition on the total wave field.
The incident plane wave 
\begin{equation}
    u^{in} (x) = \exp( ik(x_1\cos(\theta)+x_2\sin(\theta)) )
\end{equation}
impinges with the angle $\theta=135^\circ$ upon the obstacle.

\subsubsection{CG iteration and stopping criteria}
In Algorithm (Section~\ref{algo:CMCG}), the CMCG method terminates at the $\ell$-th iteration and returns 
\begin{equation}
  \ulcmcg = v_0^{(\ell)} + (i/\omega)v_1^{(\ell)} 
\end{equation}
when the \textit{relative CG-residual} in Step 5.8,
\begin{equation}
\label{eqn:stopping-criterion-cg-relres}
|\ulcmcg|_{CG}
:=
\sqrt{        
            \dfrac{
                \|\nabla r^{(\ell+1)}_{0}\|_{L^2(\Omega)}^2 
                +
                \|(1/c)\ r^{(\ell+1)}_{1}\|_{L^2(\Omega)}^2
            }
            {                  
                \|\nabla r^{(0)}_{0}\|_{L^2(\Omega)}^2 
                +
                \|(1/c)\ r^{(0)}_{1}\|_{L^2(\Omega)}^2
            }
        },
\end{equation}
is less than the tolerance $tol$. 
Indeed, a small CG-residual indicates that the gradient of $J$ is sufficiently small at $(v_0^{(\ell)},v_1^{(\ell)})$
and thus that a minimum has been reached. 

Since the cost functional $J$ also vanishes at the minimum, 
we can use $J$ itself, instead of its gradient, to monitor convergence of the CG iteration
via the \textit{relative periodicity misfit},
\begin{equation}
    \label{eqn:stopping-criterion-periodic-misfit}
    |\ulcmcg|_{J} :=     
    \frac{\sqrt{J(v_0^{(\ell)},v_1^{(\ell)})}}{\|f\|_{L^2(\Omega)}+\| g_S\|_{L^2(\Gamma_S)} }
.
\end{equation}
In fact, the convergence criterion (\ref{eqn:stopping-criterion-periodic-misfit}) is typically used in long-time simulations of the wave equation without controllability 
(``do-nothing'' approach) to determine the current misfit from periodicity in the energy norm. 

Alternatively, we may also directly compute the current \textit{relative Helmholtz residual} from (\ref{eqn:HE}):
\begin{equation}
    \label{eqn:stopping-criterion-Helmholtz-relres}
    |\ulcmcg|_{H}
    :=
    \frac{\|{\bf A}_h{{\bf u}_h^{(\ell)}}-{\bf b}_h\|_2}{\|{\bf b}_h\|_2},
\end{equation}
where ${\bf A}_h$ and ${\bf b}_h$ result from a FE discretization of (\ref{eqn:HE}) without mass-lumping,
${\bf u}_h^{(\ell)}$ corresponds to the discrete vector of FE coefficients of $\ulcmcg$, and
$\|\cdot\|_2$ denotes the discrete Euclidean norm.

\begin{figure}[t]
    \centering
    \hspace*{-1cm}
    \subfloat[
    convex obstacle
    \label{fig:square-shaped-obstacle-residual-a}
    ]{
        \includegraphics{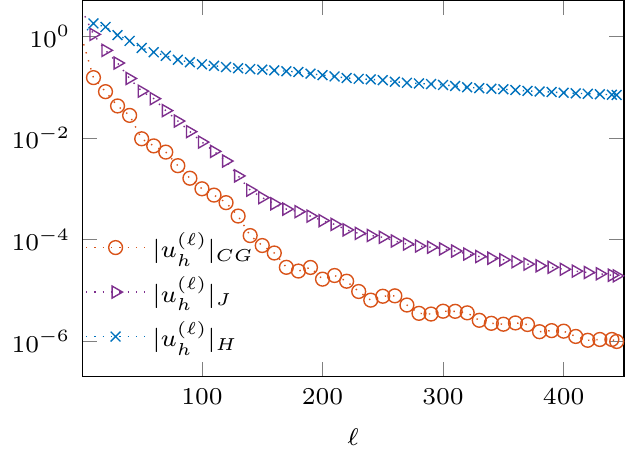}    
    }
    \hspace*{-0.5cm}
    \subfloat[nonconvex obstacle
    \label{fig:square-shaped-obstacle-residual-b}
    ]{
	\includegraphics{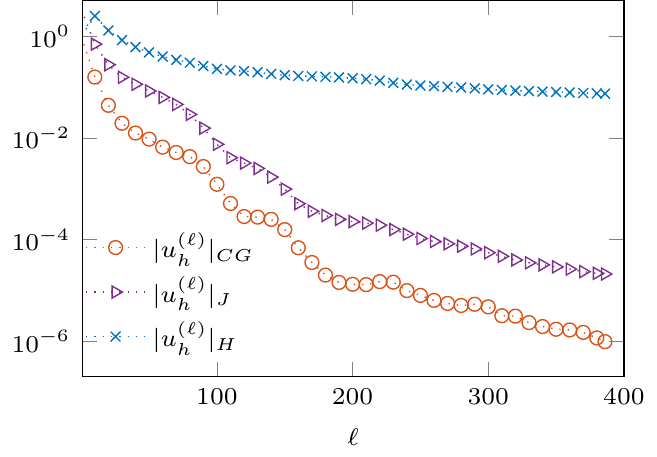}
    }
    \caption{CG iterations and stopping criteria:
    relative CG residual $|\ulcmcg|_{CG}$
    in
    (\ref{eqn:stopping-criterion-cg-relres}),
    Helmholtz residual $|\ulcmcg|_{H}$
    in
    (\ref{eqn:stopping-criterion-Helmholtz-relres}),
    and 
    periodicity mismatch $|\ulcmcg|_{J}$
    in 
    (\ref{eqn:stopping-criterion-periodic-misfit})    
    at the $\ell$-th CG iteration.
    }
    \label{fig:square-shaped-obstacle-residual}
\end{figure}

In Figure \ref{fig:square-shaped-obstacle-residual}, we monitor $|\ulcmcg|_{CG}$, $|\ulcmcg|_{J}$  and $|\ulcmcg|_{H}$,
defined in (\ref{eqn:stopping-criterion-cg-relres})--(\ref{eqn:stopping-criterion-Helmholtz-relres}) for the CMCG solution $\ulcmcg$ at the $\ell$-th CG iteration. Whether for 
a convex  (Figure \ref{fig:square-shaped-obstacle-mesh-a}) 
or  a nonconvex (Figure \ref{fig:square-shaped-obstacle-mesh-b})  obstacle, both the 
CG-residual $|\ulcmcg|_{CG}$ and the periodicity  misfit $|\ulcmcg|_{J}$ rapidly 
converge to zero.  In contrast, the Helmholtz residual $|\ulcmcg|_{H}$ stagnates beyond the first hundred CG iterations, 
as the mass-matrix that appears in  ${\bf A}_h$ in (\ref{eqn:stopping-criterion-Helmholtz-relres})  is discretized here without mass-lumping.
That additional discretization error together
with  the numerical error in the time integration of (\ref{eqn:WE})
both prevent the discrete Helmholtz residual  $|\ulcmcg|_{H}$ from converging to zero; 
hence, (\ref{eqn:stopping-criterion-Helmholtz-relres}) is generally not a reliable stopping criterion
for the CMCG method, unless the spatial FE discretizations used in (\ref{eqn:HE}) and (\ref{eqn:discrete_HE})
are identical.
\subsubsection{CMCG method vs. long-time wave equation solver}
\begin{figure}[ht]
    \centering
    \hspace*{-1cm}
    \begin{minipage}{14cm}    
    \subfloat[convex obstacle \label{fig:square-shaped-obstacle-waveequationsolver-a}]{
        \includegraphics{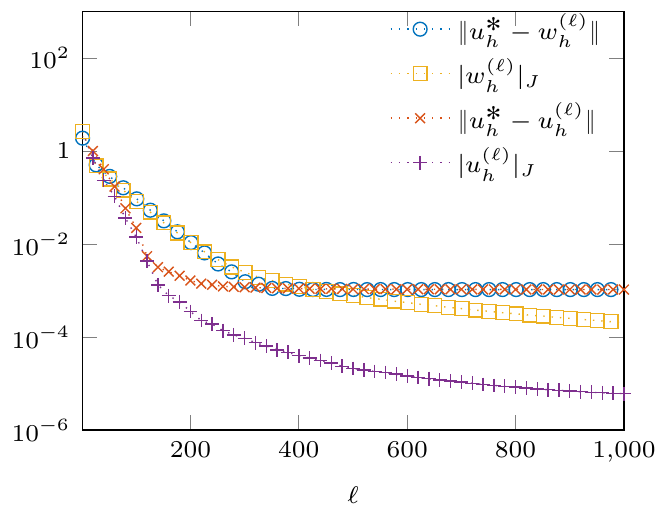}
    }
    \hspace*{-0.5cm}
    \subfloat[nonconvex obstacle \label{fig:square-shaped-obstacle-waveequationsolver-b}]{
	\includegraphics{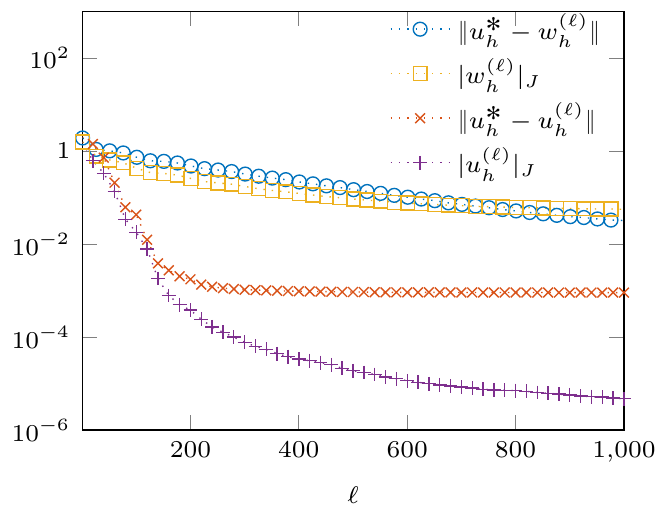}
        }
    \end{minipage}

    \caption{CMCG method vs. long-time wave equation solver: plane wave scattering from  a convex (a) or a nonconvex obstacle (b).
    Comparison between the numerical solution, $\ulcmcg$, obtained with the CMCG method at the $\ell$-th CG iteration and 
    the approximate solution $\ulasymp$,  obtained via (\ref{eqn:AP_SOLUTION}) from the solution of the wave equation at time $t=\ell\,T$
    without controllability.}    
    \label{fig:square-shaped-obstacle-waveequationsolver}
\end{figure}

In general, the solution $w(x,t)$ 
of the time-harmonically forced wave equation (\ref{eqn:WE}) converges asymptotically to the  time-harmonic solution \cite{BR1994}
\begin{equation}
    \label{eqn:AP_WAVE_SOLUTION}
    w(x,t) \sim \real{u(x)\exp(-i\omega t)} \qquad \text{as $t\rightarrow +\infty$},
\end{equation}
where $u$ is the (unique) solution of the Helmholtz equation (\ref{eqn:HE}).
Thus, with a wave equation solver at hand, one can in principle compute $u$ from $w$ by solving (\ref{eqn:WE}) without controllability
until a quasi-periodic regime is reached. 
Given the current value of $w(\cdot,t)$ at time $t=\ell\, T$, $\ell\ge 1$, 
one can extract from it the complex-valued approximate solution of (\ref{eqn:HE}),
\begin{equation}
    \label{eqn:AP_SOLUTION}
    \ulasymp
    :=
    w(\cdot,\ell T) + \frac{i}{\omega} w_t(\cdot,\ell T),
    \qquad
    \ell\ge 1,
    \qquad
    T=(2\pi)/\omega,
\end{equation}
which converges to $u$ as $\ell\rightarrow+\infty$.
This ``do-nothing'' approach only requires the time integration of (\ref{eqn:WE}) 
without controllability or CG iteration,
but it may converge  arbitrarily  slowly for nonconvex obstacles due to trapped modes \cite{BGP1998,GT2018}.

In Figure \ref{fig:square-shaped-obstacle-waveequationsolver}, we monitor 
the periodicity misfit of $|\ulcmcg|_{J}$ and 
$|\ulasymp|_{J}$, where $\ulcmcg$ is the CMCG solution at the $\ell$-th CG iteration and $\ulasymp$ is given by (\ref{eqn:AP_SOLUTION}).
In addition, we also compare both numerical solutions 
with the direct solution $\uddm$ of the linear system (\ref{eqn:discrete_HE}), 
resulting from the same underlying FE discretization, yet without mass-lumping.

We observe that the asymptotic solution $\ulasymp$
and the CMCG solution $\ulcmcg$ indeed both converge to the time-harmonic solution $\uddm$, until 
the additional errors caused by mass-lumping and the time discretization dominate the total error -- see Section \ref{sec:semi-discrete-solution}.
For the convex obstacle, the number of CG iterations required by $\ulcmcg$  is only half the number of time periods needed for $\ulasymp$  to 
reach the same level of accuracy.
However, since each CG iteration requires not only the solution of a forward and backward wave equation but also of the elliptic problem (\ref{eqn:CMCG_PRECOND_STRONGLYELLIPTIC}),
simply computing a long-time solution of the time-harmonically forced wave equation (\ref{eqn:WE}) without controllability in fact 
proves cheaper here than the CMCG Algorithm.
For a nonconvex obstacle, however, the long-time numerical solution of the time-dependent wave equation $\ulasymp$ 
converges extremely slowly and 
fails to reach the asymptotic time-harmonic regime
even after 1000 periods. In contrast, the convergence of the CMCG solution $\ulcmcg$ remains remarkably insensitive to the non-convexity of the obstacle.

\subsubsection{Initial run-up}
\begin{figure}[ht]
    \centering
    \hspace*{-2cm}
    \begin{minipage}{14cm}
    \subfloat[convex obstacle \label{fig:square-shaped-obstacle-runup-a}]{
	\includegraphics{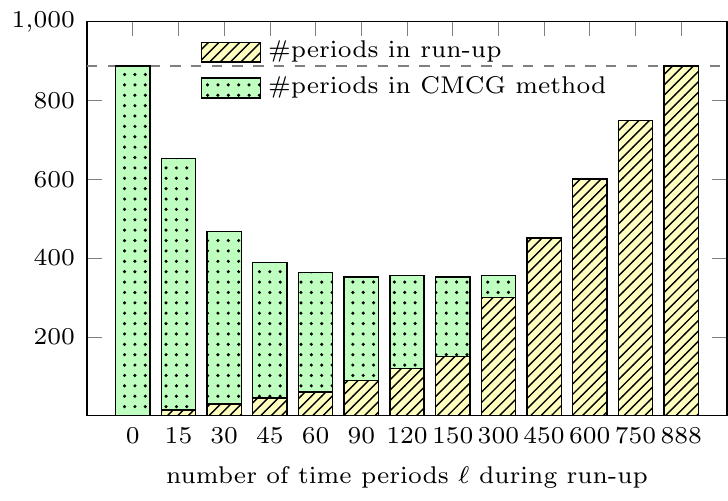}
	}
    \hspace*{-0.2cm}
    \subfloat[nonconvex obstacle \label{fig:square-shaped-obstacle-runup-b}]{
	\includegraphics{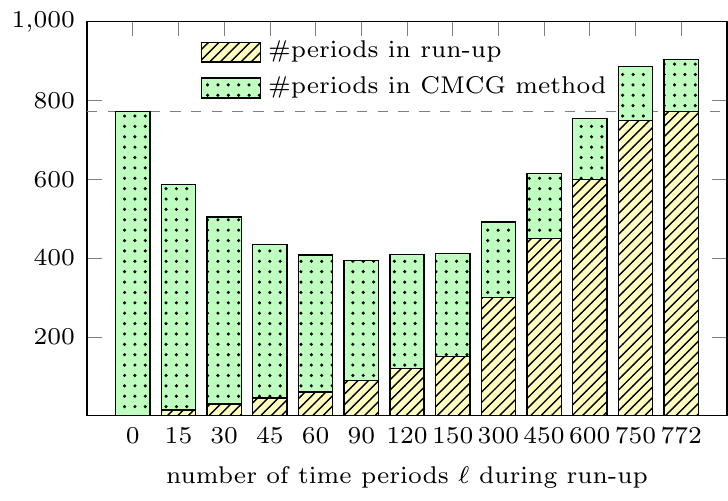}
        }
    \end{minipage}
    \caption{Initial run-up. Plane wave scattering problems from (a) a convex or (b) a nonconvex obstacle:
        total number of forward and backward wave equations solved over one period $[0,T]$ 
        until convergence.}
    \label{fig:square-shaped-obstacle-runup}
\end{figure}
In \cite{M1993}, Mur suggested that convergence of the time-harmonically forced
wave equation (\ref{eqn:WE}) to the time-harmonic asymptotic regime can be 
accelerated by pre-multiplying the time-harmonic sources in (\ref{eqn:WE}) with the smooth transient function $\theta_{tr}$ from zero to one,
\begin{equation}
    \theta_{tr}(t) = 
    \left\{
    \begin{array}{rr}
    \bigg(2-\sin\bigg(\dfrac{t}{t_{tr}}\ \dfrac{\pi}{2}\bigg)\bigg)\sin\bigg(\dfrac{t}{t_{tr}}\ \dfrac{\pi}{2}\bigg)
        ,&0\leq t \leq t_{tr},\\
     1,&t\ge t_{tr},\\
    \end{array}
    \right.
    \label{eqn:transientfunction}
\end{equation}
active during the initial time interval $[0,t_{tr}]$, 
$t_{tr}=\ell \,T$
-- see also \cite{BGP1998}.

Again, we consider plane wave scattering either from a convex  or nonconvex obstacle -- 
see Figure \ref{fig:square-shaped-obstacle-mesh}.
Now, we first solve the wave equation (\ref{eqn:WE}) with the modified source terms 
and zero initial conditions until
time $t=\ell \,T$, $\ell\ge 1$, which yields 
the time-dependent solution $y_{tr}$.
After that initial run-up phase, we then apply the CMCG Algorithm (Section~\ref{algo:CMCG}) using the initial guess 
\[
    v_0^{(0)}=y_{tr}(\cdot,\ell T), \qquad  v_1^{(0)}=(y_{tr})_{t}(\cdot,\ell T).
\]
To estimate the total computational effort,
we count the total number of time periods 
for which the (forward or backward) 
wave equation is solved: $\ell$  during initial run-up 
and $2\times\#iter_{CG}$ during the CG iteration.
In Figure \ref{fig:square-shaped-obstacle-runup}
we display the total number $2\times\#iter_{CG}+\ell$ of time periods
needed until convergence with $tol =10^{-6}$, as we vary
 the number of periods $\ell$ in the initial run-up.

For a convex obstacle, the CMCG Algorithm without any initial
run-up requires $888$ time periods. However, as in Section
\ref{sec:square-shaped-scattering-problem}, convergence can
also be achieved at a comparable computational effort
simply by solving the wave equation, here with 
the source terms pre-multiplied by $\theta_{tr}$ in (\ref{eqn:transientfunction}).
Still, the minimal computational cost is achieved when both the initial 
run-up and the CMCG Algorithm are combined.

For the nonconvex obstacle, however, simply solving the time-harmonically forced 
wave equation over a very long time, be it with or without $\theta_{tr}(t)$ smoothing,
fails to reach the long-time asymptotic final time-harmonic state.
Regardless of the length of the initial run-up, convergence indeed cannot be achieved here (within
1000 time periods) without controllability because of trapped modes. Nevertheless,
the initial run-up always speeds up the convergence
of the CMCG method by providing a judicious initial guess for the CG iteration.

\subsection{Parallel computations\label{sec:DDM}}
Both the CMCG method for the second-order formulation from Section 2 and that for the first-order formulation from Section 3 lead to inherently parallel non-intrusive algorithms, as long as an efficient parallel solver for the  time-dependent wave equation is available. As the first-order formulation with the HDG discretization neither requires mass-lumping nor the solution of an elliptic problem, it is in fact trivially parallel. 
Here we demonstrate that even the CMCG approach for the second-order formulation, which does require
the solution of  (\ref{eqn:CMCG_PRECOND_STRONGLYELLIPTIC}) at each CG iteration, nonetheless achieves strong scalability on a massively parallel architecture.

The CMCG  Algorithm from Section~\ref{algo:CMCG}  is implemented within FreeFem++~\cite{Hecht:2012:NDF}, 
an open source finite element software written in C++. FreeFem++ defines a high-level Domain Specific Language (DSL) and natively supports distributed parallelism with MPI. 
The parallel implementation of the
CMCG method relies on the spatial decomposition of the computational domain $\Omega $ into
multiple subdomains, each assigned to a single computing core. 
Local finite element spaces are then defined on the local meshes of the subdomains, 
effectively distributing the global set of degrees of freedom across the available cores. 

The bulk of the computational work for solving the forward and backward wave equations in Step~\ref{cmcgalgo:wave} of the CMCG Algorithm simply
consists in performing a sparse matrix-vector product at each time step, 
which is easily parallelized in this domain decomposition framework: it amounts to performing local matrix-vector products in parallel on the 
local set of degrees of freedom corresponding to each subdomain, 
followed by local exchange of shared values between neighboring subdomains. 

While the explicit time integration of the wave equation is trivially parallelized thanks to mass-lumping,
achieving good parallel scalability for the elliptic problem in Step~\ref{cmcgalgo:elliptic} of the CMCG Algorithm is more difficult. Here we use domain decomposition (DD) methods \cite{DJN2015}, 
which are well-known
to produce robust and scalable parallel preconditioners for the iterative solution of large scale partial differential equations. 
We use the parallel DD library HPDDM~\cite{jolivet2013scalable}, 
which implements efficiently various 
Schwarz and substructuring methods 
in C++11 with MPI and OpenMP for parallelism and 
is interfaced with FreeFem++
.

The elliptic problem (\ref{eqn:CMCG_PRECOND_STRONGLYELLIPTIC})  in the CMCG algorithm is solved by HPDDM using a two-level overlapping Schwarz DD preconditioner, 
where the coarse space is built using Generalized Eigenproblems in the Overlap (GenEO)~\cite{Spillane:2014:ASC}. 
The GenEO approach has proved effective in producing highly scalable preconditioners for solving various elliptic problems ~\cite{BDGST2017,Spillane:2014:ASC}. 

All computations were performed on the supercomputer \textit{O{\small CCIGEN}} at \textit{C\small{INES}}, France~\footnote{\url{https://www.cines.fr/calcul/materiels/occigen/}},
with 50544 (\textit{Intel XEON Haswell}) cores.

\subsubsection{2D Marmousi Model}
\begin{figure}[t]
    \centering
    \includegraphics[width=12cm]{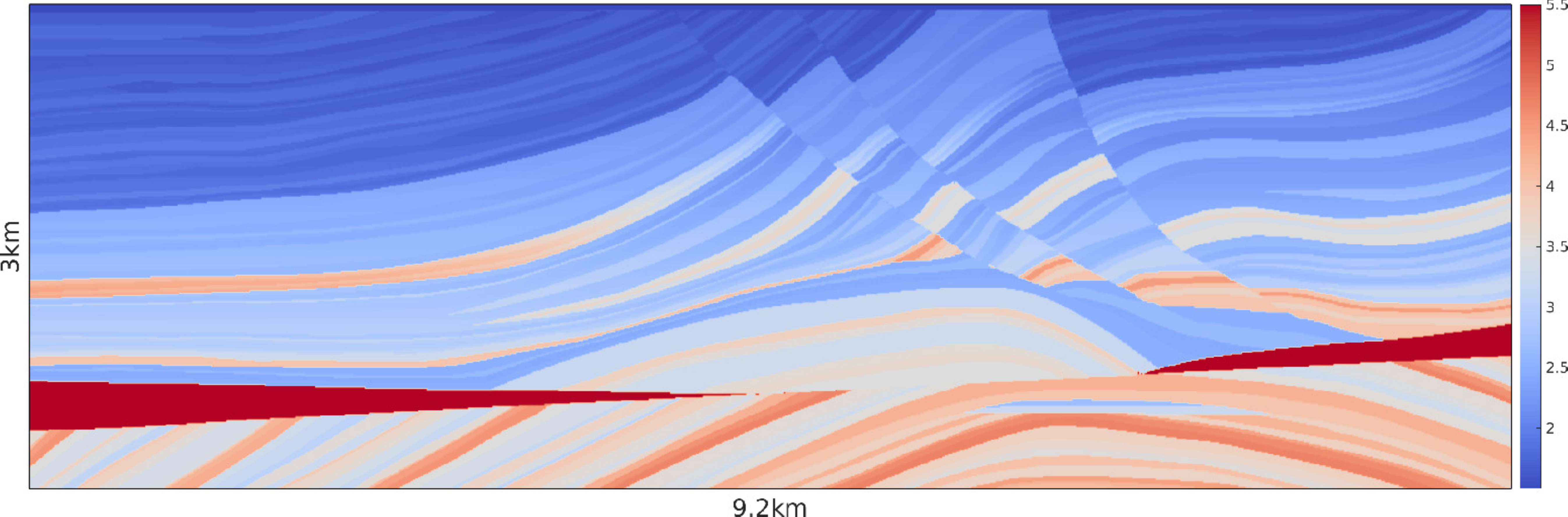} 
    \caption{Marmousi model: propagation velocity $1.5\leq c(x)\leq 5.5$  $[km/s]$}
    \label{fig:marmousi-profile}
\end{figure}
\begin{table}[t]
\centering
\small
\hspace*{-0.5cm}
\begin{tabular}{r|r|r|r}
{\bf Frequency} & {\bf Wave number}  &{\bf \#Unknowns}  & {\bf\#Nodes}\\
$\nu$ [Hz] & $k=\omega/c=2\pi \nu/c$&$ndof$            & 24 cores per node\\\hline
$10$  & $11$ -- $42$    & $    1'658'443$ &  $1$--$8$\\
$20$  & $22$ -- $84$    & $    6'628'881$ & $1$--$16$\\
$40$  & $45$ -- $168$   & $   26'505'761$ & $8$--$64$\\
$60$  & $68$ -- $252$   & $   59'630'641$ & $16$--$128$\\
$80$  & $91$ -- $336$   & $  106'003'521$ & $16$--$128$\\
$160$ & $182$ -- $671$  & $  423'975'041$ & $64$--$256$\\
$250$ & $285$ -- $1048$ & $1'035'241'009$ & $128$--$512$\\
\end{tabular}
\caption{2D-Marmousi model: $P^2$-FE
    with $15$ points per wave length}
\label{tbl:2d-marmousi-configuration}
\end{table}
\begin{figure}[t]
    \centering    
    \includegraphics[width=12cm,height=5.55cm]{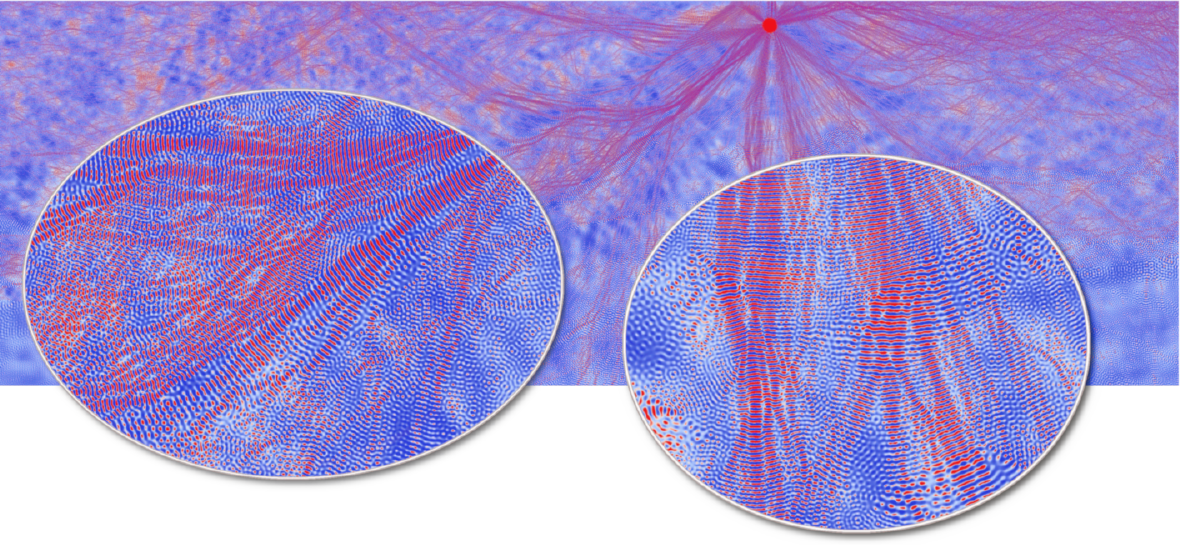}
    \caption{2D-Marmousi model.
        Real part of the wave field        
        with $\omega=2\pi \nu$, $\nu=250$ [Hz]       
    }    
    \label{fig:2d-marmousi-solution}
\end{figure}
\begin{figure}[t]        
    \centering
    \hspace*{-2cm}
    \includegraphics{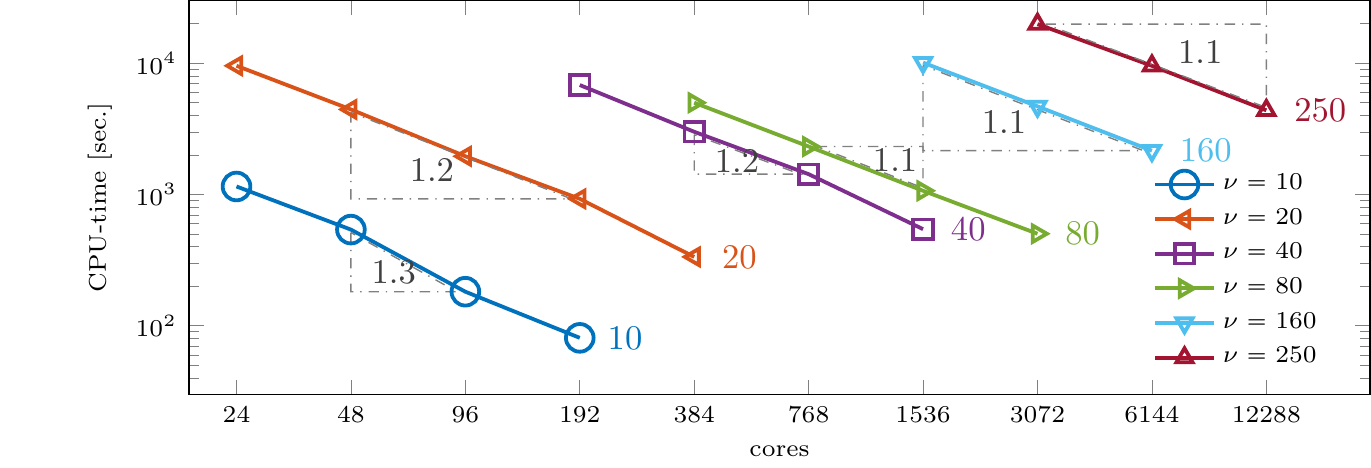}
    \caption{2D-Marmousi model.
    Total CPU-time in seconds 
    for varying number of cores.
    For each frequency $\nu$,
    the FE-discretization 
    and problem size remain fixed.
    }    
    \label{fig:2d-marmousi-cputtime}
\end{figure}

Here we consider the well-known Marmousi model from geophysics \cite{BBLPP1990},
that is (\ref{eqn:HE}) in $\Omega = (0,9.2) \times (0,3)$ $[km]$ with the source
\[
    f(x) =  \exp(-2000((x-x_1^*)^2+(y-x_2^*)^2)),
    \qquad
    (x_1^*,x_2^*) = (6, -3/16).
\]
The velocity profile $c(x)$ is shown in Figure \ref{fig:marmousi-profile}
and we apply absorbing boundary conditions on the lateral and lower boundaries
and a homogeneous Dirichlet condition at the top.
For the spatial discretization, we use a $P^2$-FE method with 
(order preserving) mass-lumping \cite{CJRT2001} and at least $15$ points per wave length.
For the time integration of (\ref{eqn:WE}), we apply the leap-frog scheme (LF); here,
the number of $T/\Delta t=390$ time steps per period remains constant at all frequencies $\nu = \omega/2\pi$, 
as both $T$ and $\Delta t$ are inversely proportional to $\nu$.
To speed-up the convergence of the CMCG method, we also use an initial run-up (Section \ref{sec:square-shaped-scattering-problem}) until time $t_{tr}$, which lets
waves travel at least once across the entire computational domain
during run-up; hence, we set
\[
    \ell = 
    \left\lceil
    \frac{\sqrt{9.2^2+3^2}}{Tc_{\min}}    
    \right\rceil,
    \qquad
    t_{tr}=\ell \,T,
    \qquad
    T=(2\pi)/\omega.
\]

For any particular frequency $\nu$,
we apply the CMCG method for fixed parameters and FE-mesh while increasing the number of (CPU) cores.
Figure \ref{fig:2d-marmousi-solution} displays the real part of the wave field with $\nu=250$ [Hz].
In Figure \ref{fig:2d-marmousi-cputtime},
we observe linear speed-up (strong scaling) at every frequency
with increasing number of cores. 
In fact, the speed-up is even slightly better than linear due to cache effects, but also because
the cost of the direct solver used on each subdomain decreases superlinearly with the decreasing 
size of subdomains as the number of cores increases.

As the  frequency $\nu$ increases, both the period $T=1/\nu$ 
and the time-step $\Delta t$ decrease, so that the number of time steps per CG iteration remains constant.
Since the number of CG iterations does not grow here with increasing $\nu$,
the bulk of the computational work in the CMCG Algorithm in fact shifts to the run-up phase. 
For $\nu=10$ Hz, for instance, the CMCG Algorithm stops after $273$ CG iterations,
while $74\%$ of the total computational time is spent in the time integration of (\ref{eqn:WE}), 
$16\%$ in the elliptic solver (DDM) and $10\%$ in the initial run-up.
In contrast, for $\nu=250$ Hz, the CMCG Algorithm already stops after $5$ CG iterations,
while $99\%$ of the total computational time is spent in the initial run-up and $1\%$ in the CG iteration.
By modifying the run-up time $t_{tr}$, one could arbitrarily shift the relative computational cost between 
run-up and CG iterations and thus further optimize for a minimal total execution time.
 \subsubsection{3D cavity} 
\begin{figure}[t]
    \centering   
    \hspace*{-0.25cm}
    \subfloat[front view]{        
        \includegraphics{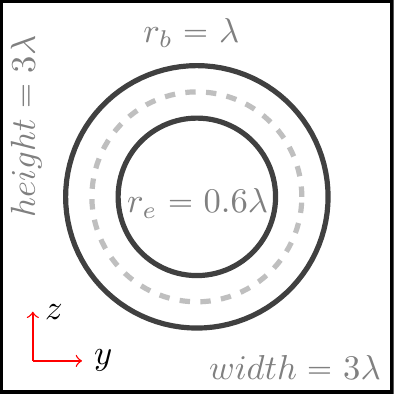}
    }
    \hspace*{-0.25cm}
    \subfloat[cross-section]{
        \includegraphics{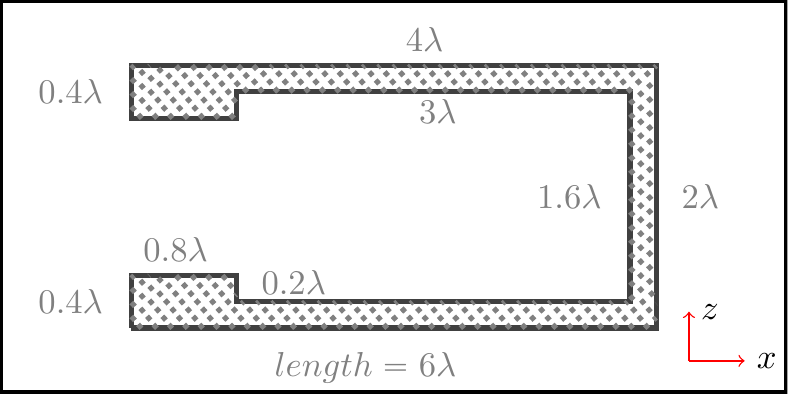}
    }
    \caption{3D-cavity:
    a) front view of the opening with inner and outer radius,
    b) longitudinal cross-section.
    }
    \label{fig:3d-scattering-mesh}
\end{figure}    
\begin{figure}[t]    
\centering
    \includegraphics[width=10cm]{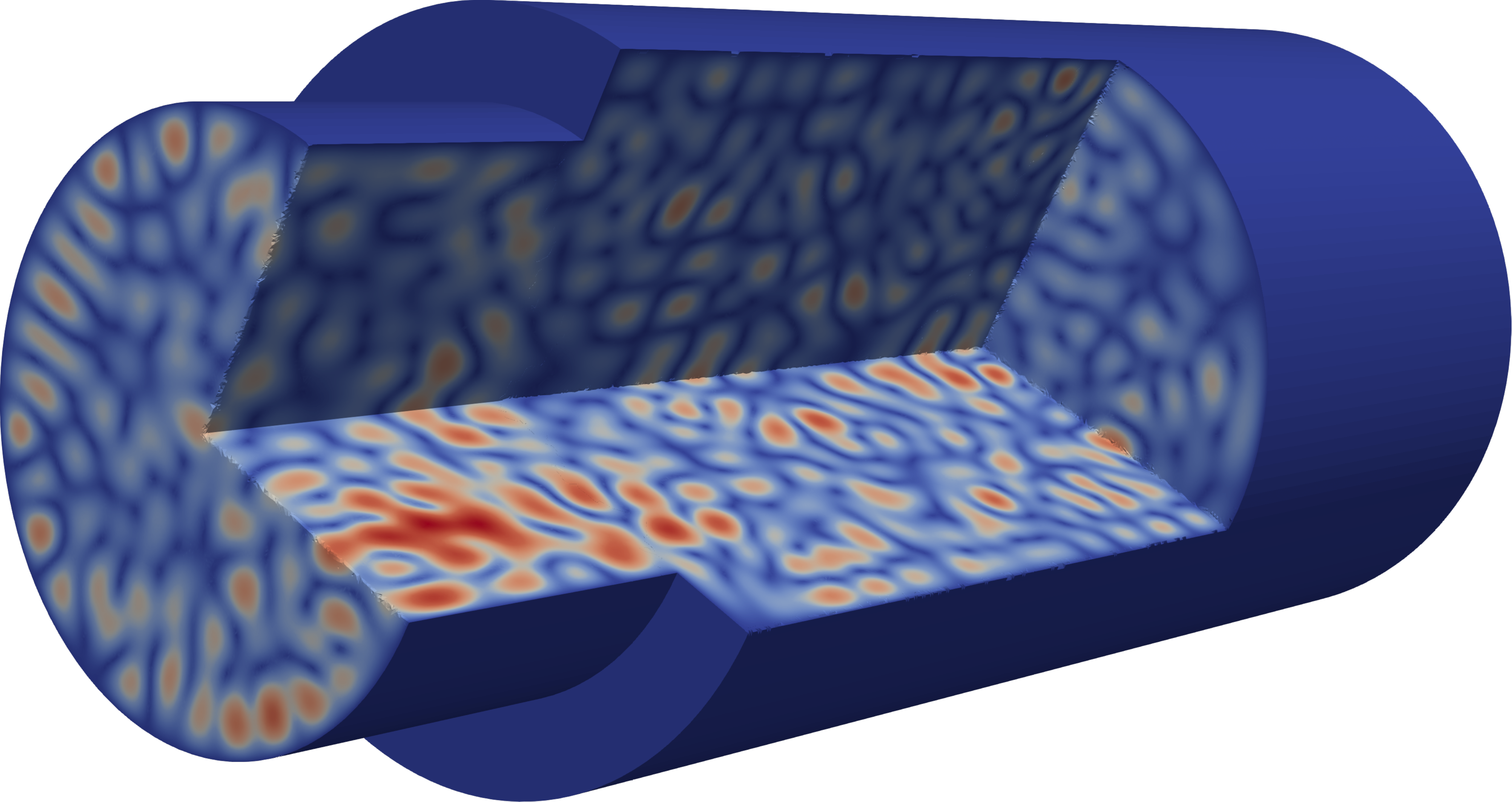}    
    \caption{3D-cavity; total wave field (\ref{eqn:HE}) with $c=1$, $\omega=2\pi \nu$ and $\nu=6$ obtained with the CMCG method
    }
    \label{fig:3d-scattering-solution}
\end{figure}
\begin{table}[t]
\centering
\small  
\hspace*{-1cm}
\begin{tabular}{c|r|r|r|r}
{\bf Frequency} & {\bf\#Unknowns}  & {\bf\#Tetrahedra} &   {\bf CG iterations}  & {\bf\#Nodes} 
\\
$\nu=2\pi \omega$   & $ndof$          &                              & &24 cores per node
\\\hline
$2$ & $8.17\cdot  10^5$  & $5'051'049$    & $239$  & $1$--$8$\\
$3$ & $5.22\cdot  10^6$  & $31'190'000$   & $440$  & $2$--$32$\\
$4$ & $1.9 \cdot  10^7$  & $114'391'112$  & $607$  & $32$--$96$\\
$6$ & $1.18\cdot  10^8$  & $703'590'464$  & $578$  & $64$--$128$\\
\end{tabular}
\caption{3D-cavity: CMCG methods with $P^1$-FEM.
As $\eta$ increases, the ratio
$h k^{3/2}$  remains constant to avoid pollution errors \cite{BS1997}.}
\label{tbl:3d-scattering-configuration}
\end{table}
\begin{figure}[t]
    \centering   
    \hspace*{-2cm}
    \includegraphics{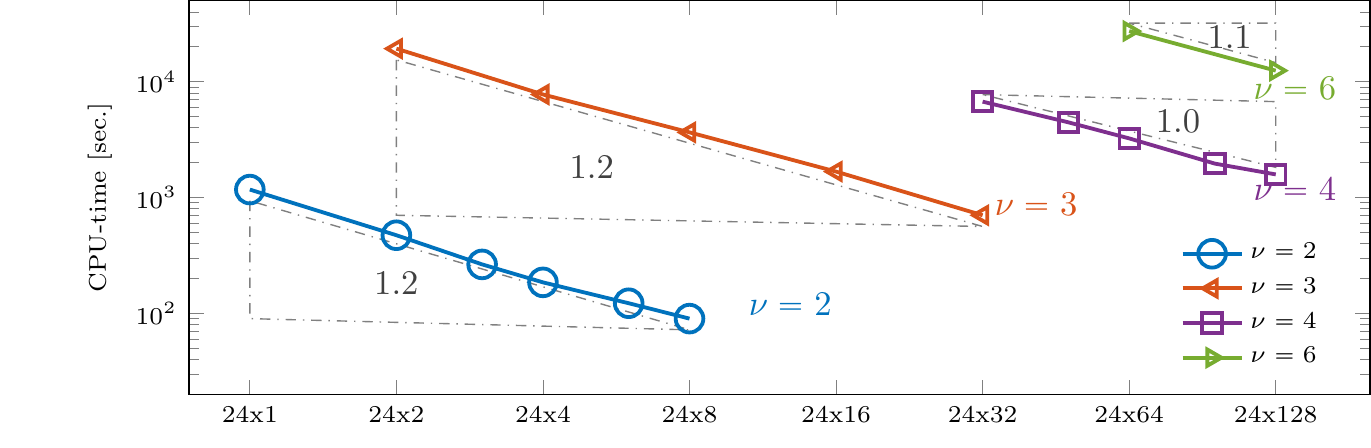}
    \caption{3D-cavity: Total CPU-time in seconds 
    for varying number of cores.
     For each frequency $\nu$,
    the FE-discretization 
    and problem size remain fixed.}
    \label{fig:3d-scattering-timing}
\end{figure}
Finally, we compute the scattered wave from a sound-soft cavity -- see Figure \ref{fig:3d-scattering-mesh} -- 
and hence consider (\ref{eqn:HE}) in $\Omega = (0,6) \times (0,3)\times(0,3)$ with $c=1$, $k=\omega=2\pi\nu$, $\lambda=1$,
$f\equiv g_D\equiv g_N \equiv 0$ and
\[
    g_S=-(\partial_n-ik)u^{in}, \qquad  u^{in}({\bf x}) = \exp( ik\ {\bf x}^\intercal {\bf d}) ,
\qquad {\bf d}=(1/2,0,\sqrt3/2)^\intercal
.
\]
We impose a homogeneous Dirichlet boundary condition on the obstacle 
and a Sommerfeld-like absorbing condition on the exterior boundary of $\Omega$.

Now, we discretize (\ref{eqn:WE}) with $P^1$-FE in space
and the second-order LF method in time.
To control the pollution error, we set $hk^{3/2} \sim const$,
as we increase the frequency $\nu$.
Figure \ref{fig:3d-scattering-solution} shows the total wave field with $\nu=6$ inside the cavity.
For fixed parameters and mesh size, we now solve (\ref{eqn:HE}) at 
frequencies $\nu=2,3,4,6$
with the CMCG method using an increasing number of cores -- see Table \ref{tbl:3d-scattering-configuration}.
Again, we observe  in Figure \ref{fig:3d-scattering-timing} (better than) linear (strong) scaling with increasing number  of cores. In contrast to the previous Marmousi problem, the
"do-nothing" approach without controllability fails here because the 3D cavity is not convex.

\section{Concluding remarks}
We have presented two inherently parallel controllability methods (CM) for the numerical solution of the
Helmholtz equation in heterogeneous media. The first, based on the
second-order formulation of the wave equation, uses a standard (continuous) FE discretization in space with order preserving mass-lumping. Each conjugate gradient (CG) iteration then
requires the explicit time integration of a forward and backward wave equation, together with
the solution of the symmetric and coercive elliptic problem (\ref{eqn:CMCG_PRECONDITIONIER}), which is independent of the frequency.
The second, based on the first-order (or mixed) formulation of the wave equation, uses a recent
hybridized discontinuous Galerkin (HDG) discretization, which not only automatically yields a block-diagonal mass-matrix
but also completely avoids solving (\ref{eqn:CMCG_PRECONDITIONIER}). Hence, it is trivially parallelized and even leads to superconvergence after a local post-processing step.

Both CMCG methods are inherently parallel, as they lead to iterative algorithms whose
convergence rate is independent of the number of cores on a distributed memory architecture.
Thanks to the well-known parallel efficiency of explicit methods
combined with the excellent scalability of 
two-level domain decomposition preconditioners for coercive 
elliptic problems up to thousands of cores implemented in HPDDM, even the second-order CMCG approach exhibits parallel strong scalability.

The CMCG method can be applied to general boundary-value problems governed by the Helmholtz equation, such as sound-soft or sound-hard scattering problems or wave propagation in physically bounded domains. Although the CMCG solution will generally contain higher order spurious eigenmodes, 
we have proposed in Section \ref{sec:FilterFundamentalFrequency} a simple filtering procedure to remove them.
Furthermore, including a transient initial run-up to determine a judicious initial guess
significantly accelerates the CG iteration. In fact, for scattering from convex obstacles, simply solving the time-harmonically forced wave equation over a long-time without any controllability
can provide an even simpler, highly parallel Helmholtz solver.
For nonconvex obstacles, however, solving the wave equation without any controllability  ("do-nothing" approach)
is not a viable option, as the long time asymptotic convergence to the time-harmonic regime is simply too slow due to trapped modes. In all cases, 
the CMCG Algorithm combined with the initial run-up leads
to the smallest time-to-solution.

The CMCG approach developed here for the
Helmholtz equation immediately generalizes to other time-harmonic vector wave equations from
electromagnetics or elasticity. Its implementation is non-intrusive and 
particularly useful when a parallel efficient time-dependent wave equation solver is at hand.
In the presence of local mesh refinement, local time-stepping
methods \cite{GPRS2017} permit to circumvent the increasingly stringent CFL condition without sacrificing the explicitness or inherent parallelism. 
Finally, the CMCG method can also be used to compute periodic, but not necessarily time-harmonic, 
solutions of the wave equations.
In particular, if the source consists of a superposition of several time-harmonic sources ("super-shot") with rational frequencies, the solutions to the different Helmholtz problems can be extracted via filtering from a 
single application of the CMCG method.

\smallskip
{\bf Acknowledgement:}  This work was supported by the Swiss National Science Foundation
 under grant SNF 200021\_169243. Access to the HPC resources of CINES was granted under allocation 2018-A0040607330 by GENCI.

\bibliographystyle{elsarticle-num}
\bibliography{HPCControllabilityMethods}

\begin{thebibliography}{10}
\expandafter\ifx\csname url\endcsname\relax
  \def\url#1{\texttt{#1}}\fi
\expandafter\ifx\csname urlprefix\endcsname\relax\def\urlprefix{URL }\fi
\expandafter\ifx\csname href\endcsname\relax
  \def\href#1#2{#2} \def\path#1{#1}\fi

\bibitem{EG2005}
O.~G. Ernst, M.~J. Gander, Why it is Difficult to Solve {H}elmholtz Problems
  with Classical Iterative Methods, Springer Berlin Heidelberg, Berlin,
  Heidelberg, 2012, pp. 325--363.

\bibitem{ERLANGGA2004409}
Y.~Erlangga, C.~Vuik, C.~Oosterlee, On a class of preconditioners for solving
  the {H}elmholtz equation, Appl. Num. Mat. 50~(3) (2004) 409--425.

\bibitem{Calandra:2017:GMP}
H.~Calandra, S.~Gratton, X.~Vasseur, A Geometric Multigrid Preconditioner for
  the Solution of the {H}elmholtz Equation in Three-Dimensional Heterogeneous
  Media on Massively Parallel Computers, Springer Internat. Publ., 2017, pp.
  141--155.

\bibitem{BGS2009}
M.~Bollh{\" o}fer, M.~J. Grote, O.~Schenk, Algebraic multilevel preconditioner
  for the {H}elmholtz equation in heterogeneous media, SIAM J. Sci. Comput.
  31~(5) (2009) 3781--3805.

\bibitem{graham2017domain}
I.~Graham, E.~Spence, E.~Vainikko, Domain decomposition preconditioning for
  high-frequency {H}elmholtz problems with absorption, Mathematics of
  Computation 86~(307) (2017) 2089--2127.

\bibitem{BDGST2017}
M.~Bonazzoli, V.~Dolean, I.~G. Graham, E.~A. Spence, P.-H. Tournier, A
  two-level domain-decomposition preconditioner for the time-harmonic
  {M}axwell's equations, Lect. Notes Comput. Sci. Eng.

\bibitem{EY2011}
B.~Engquist, L.~Ying, Sweeping preconditioner for the {H}elmholtz equation:
  Moving perfectly matched layers, Mult. Model. Sim. 9 (2011) 686--710.

\bibitem{BGP1998}
M.-O. Bristeau, R.~Glowinski, J.~P{\'e}riaux, Controllability {M}ethods for the
  {C}alculation of {T}ime-{P}eriodic {S}olutions. {A}pplication to
  {S}cattering, J. Comput. Phys. 147~(2) (1998) 265--292.

\bibitem{HMPR2007}
E.~Heikkola, S.~M{\"o}nk{\"o}l{\"a}, A.~Pennanen, T.~Rossi, Controllability
  method for acoustic scattering with spectral elements, J. Comput. Appl. Math.
  204~(2) (2007) 344--355.

\bibitem{HMPR2007_2}
E.~Heikkola, S.~M{\"o}nk{\"o}l{\"a}, A.~Pennanen, T.~Rossi, Controllability
  method for the {H}elmholtz equation with higher-order discretizations, J.
  Comput. Phys. 225~(2) (2007) 1553--1576.

\bibitem{GT2018}
M.~J. Grote, J.~H. Tang, On controllability methods for the {H}elmholtz
  equation, tech. report 2018-06, University of Basel, 2018 (2018).

\bibitem{GR2006}
R.~Glowinski, T.~Rossi, A mixed formulation and exact controllability approach
  for the computation of the periodic solutions of the scalar wave
  equation.(i): Controllability problem formulation and related iterative
  solution., Comptes Rendus Math. 343~(7) (2006) 493--498.

\bibitem{CNPS2016}
B.~Cockburn, N.~Nguyen, J.~Peraire, M.~Stanglmeier, An explicit hybridizable
  discontinuous {G}alerkin method for the acoustic wave equation, Comput. Meth.
  Appl. Mech. Engrg. 300 (2016) 748--769.

\bibitem{L1988}
J.-L. Lions, Exact controllability, stabilization and perturbations for
  distributed systems, SIAM J. Appl. Math. 30~(2) (1988) 1--68.

\bibitem{BGT1982}
A.~Bayliss, M.~Gunzburger, E.~Turkel, Boundary conditions for the numerical
  solution of elliptic equations in exterior region, SIAM J. Appl. Math. 42~(2)
  (1982) 430--451.

\bibitem{GK1995}
M.~J. Grote, J.~B. Keller, {O}n nonreflecting boundary conditions, J. Comput.
  Phys. 122~(2) (1995) 231--243.

\bibitem{MS2014}
J.~M\'{a}lek, Z.~Strako\v{s}, Preconditioning and the Conjugate Gradient Method
  in the Context of Solving {PDE}s, SIAM, 2014.

\bibitem{DJN2015}
V.~Dolean, P.~Jolivet, F.~Nataf, An Introduction to Domain Decomposition
  Methods. Algorithms, Theory, and Parallel Implementation, SIAM, 2015.

\bibitem{E2010}
L.~C. Evans, Partial Differential Equation, AMS, 2010.

\bibitem{P1983}
A.~Pazy, Semigroups of Linear Operators and Applications to Partial
  Differential Equations, Springer, 1983.

\bibitem{CX2006}
P.~Cummings, X.~Feng, Sharp regularity coefficient estimates for complex-valued
  acoustic and elastic {H}elmholtz equations, Math. Models Methods Appl. Sci.
  16~(01) (2006) 139--160.

\bibitem{BJT2000}
E.~B{\'e}cache, P.~Joly, C.~Tsogka, An analysis of new mixed finite elements
  for the approximation of wave propagation problems, SIAM J. on Numer. Anal.
  37~(4) (2000) 1053--1084.

\bibitem{CJRT2001}
G.~Cohen, P.~Joly, J.~E. Roberts, N.~Tordjman, Higher order triangular finite
  elements with mass lumping for the wave equation, SIAM J. Numer. Anal. 38~(6)
  (2001) 2047--2078.

\bibitem{BR1994}
C.~Bardos, J.~Rauch, Variational algorithms for {H}elmholtz equation using time
  evolution and artificial boundaries, Asympt. Anal. 9 (1994) 101--117.

\bibitem{M1993}
G.~Mur, The finite-element modeling of three-dimensional electromagnetic fields
  using edge and nodal elements, IEEE Trans. on Antenn. and Prop. 41 (1993)
  948--953.

\bibitem{Hecht:2012:NDF}
F.~Hecht, New development in {F}ree{F}em++, J. Num. Math. 20 (2012) 251--265.

\bibitem{jolivet2013scalable}
P.~Jolivet, F.~Hecht, F.~Nataf, C.~Prud'homme, Scalable domain decomposition
  preconditioners for heterogeneous elliptic problems, in: Proc. 2013 ACM/IEEE
  Conf. on Supercomputing, SC13, ACM, 2013, pp. 1--11.

\bibitem{Spillane:2014:ASC}
N.~Spillane, V.~Dolean, P.~Hauret, F.~Nataf, C.~Pechstein, R.~Scheichl,
  Abstract robust coarse spaces for systems of {PDE}s via generalized
  eigenproblems in the overlaps, Numer. Math. 126~(4) (2014) 741--770.

\bibitem{BBLPP1990}
A.~Bourgeois, M.~Bourget, P.~Lailly, M.~Poulet, P.~Ricarte, R.~Versteeg,
  Marmousi, model and data, 1990.

\bibitem{BS1997}
I.~M. Babu{\v s}ka, S.~A. Sauter, Is the pollution effect of the {FEM}
  avoidable for the {H}elmholtz equation considering high wave numbers, SIAM J.
  Numer. Anal. 34~(6) (1997) 2392--2423.

\bibitem{GPRS2017}
M.~J. Grote, D.~Peter, M.~Rietmann, O.~Schenk, Newmark local time stepping on
  high-performance computing architectures, J. Comput. Phys. 334 (2017)
  308--326.

\end{thebibliography}
\end{document}